\documentclass[12pt,a4paper]{article}
\usepackage{amsmath,amssymb,latexsym}
\usepackage{amsthm}
\usepackage{bm}
\usepackage{overpic}
\usepackage{color}
\usepackage{xcolor}
\numberwithin{equation}{section}

\theoremstyle{plain}% default
\newtheorem{theorem}{Theorem}[section]

\newtheorem{lemma}[theorem]{Lemma}
\theoremstyle{definition}

\theoremstyle{remark}
\newtheorem{remark}[theorem]{Remark}

\def\XXint#1#2#3{{\setbox0=\hbox{$#1{#2#3}{\int}$}
\vcenter{\hbox{$#2#3$}}\kern-.5\wd0}}

\title{The error bounds of Gauss quadrature formulae for the modified weight functions of Chebyshev type}

\author{Ram\'{o}n Orive, Aleksandar V. Pej\v cev, Miodrag M. Spalevi\'c}

\date{\today}

\begin{document}

\vspace{1cm} \maketitle

%\begin{document}
%\begin{frontmatter}

\author{\bigskip\em Dedicated to Prof. Gradimir V. Milovanovi\' c on the occasion of his 70-th birthday}

\begin{abstract}In this paper, we consider the Gauss quadrature formulae corresponding to some modifications of anyone of the four Chebyshev weights, considered by Gautschi and Li in \cite{gauli}. As it is well known, in the case of analytic integrands, the error of these
	quadrature formulas can be represented as a contour integral with a
	complex kernel. We study the kernel, as it is often considered, on elliptic contours with foci
	at the points $\mp 1$ and such that the sum of semi-axes is $\rho>1$, of the
	mentioned quadrature formulas, and derive some error bounds for them. In addition, we obtain, for the first time as far as we know, a result about the behavior of the modulus of the corresponding kernels on those ellipses in some cases.
	Numerical examples checking the accuracy of such error bounds are included.
\end{abstract}

\textbf{keywords.} Gauss quadrature formulae, Chebyshev weight functions, contour integral representation,
		remainder term for analytic
		functions, error bound

\medskip

%\end{keyword}

%\end{frontmatter}

\section{Introduction}
\setcounter{equation}{0} \setcounter{theorem}{0}
\setcounter{figure}{0} \setcounter{table}{0}

In \cite{gauli}, the authors studied orthogonal polynomials with respect to some particular polynomial modifications of a given measure. Namely, given a positive measure $d\sigma$ on the real axis, with $\displaystyle \{\pi_n\}$ being its corresponding sequence of orthogonal polynomials, they considered the modified measure $\displaystyle d\widehat{\sigma}_n = \pi_n^2\,d\sigma\,,$ $n$ being an arbitrarily fixed  nonnegative integer, referring to the related orthogonal polynomials $\displaystyle \{\widehat{\pi}_{m,n}\}$ as ``induced'' orthogonal polynomials. As pointed out by the authors in \cite{gauli}, this kind of modifications of measures take place, for instance, when dealing with constrained polynomial least squares approximation (see e.g. \cite{gau89}), or in terms of providing additional interpolation points (the zeros of the induced polynomial $\displaystyle \{\widehat{\pi}_{n+1,n}\}$) in the process of extending Lagrange interpolation at the zeros of $\pi_n$ (see \cite{Bellen}).

Taking into account these and other applications, it seems natural to consider the numerical computation of integrals of the form
\begin{equation*}\label{integral}
I_{\sigma} (f) = I(f;\sigma,n) = \int f(t)\,d\widehat{\sigma}_n (t)
\end{equation*}
by means of quadrature formulae; in particular, Gauss quadratures are our main subject of interest. It is well known that the zeros and nodes of the Gauss rule can be efficiently computed by means of the eigenvalues and eigenvectors of the related tridiagonal Jacobi matrix. In general, it is not feasible to get closed analytic expressions of the entries of the Jacobi matrix for the induced measure $\displaystyle d\widehat{\sigma}_n$ in terms of the corresponding for $d\sigma$; in this sense, in \cite{gauli} a stable numerical algorithm is given. But in the particular case of the well--known four Chebyshev weights $\displaystyle d\sigma^{[i]}\,,\,i=1,2,3,4\,,$ the related induced orthogonal polynomials $\displaystyle \{\widehat{\pi}_{m,n}^{[i]}\}$ are easily expressible as combinations of Chebyshev polynomials of the first kind $T_k$ (i.e., orthogonal polynomials with respect to the Chebyshev weight $d\sigma^{[1]}\,$ (see \cite[\S 3]{gauli}). These last results will be very useful for our analysis of the error of the related quadrature formulas.

Our analysis of the error is based on its well known representation in terms of an integral contour of an appropriate kernel.  Namely, if we use a Gaussian rule $I_m(f)\,,$ with $m$ nodes, to approximate the value of the integral $I_w (f)$, for a certain positive weight function on a compact real interval, say $[-1,1]$, and an analytic integrand $f$ in a neighborhood $\Omega$ of this interval, the error of quadrature admits the following integral representation
\begin{equation}\label{intrepres}
R_m(f) = I_w (f) - I_m(f) = \,\frac{1}{2\pi i}\,\oint_{\Gamma}\,K_m(z)\,f(z)\,dz\,,
\end{equation}
where the kernel $K_m$ is given by
\begin{equation}\label{kernel}
K_m (z) = \,\frac{\varrho_m (z)}{\pi_m (z)}\,,\;\;\varrho_m (z) = \int_{-1}^1\,\frac{\pi_m (t)}{z-t}\,w(t)\,dt\,,
\end{equation}
with $\pi_m$ denoting, as usual, the $m$-th orthogonal polynomial with respect to $w$ and $\Gamma \subset \Omega$ being any closed smooth contour surrounding the real interval $[-1,1]$. Since they are the level curves for the conformal function which maps the exterior of $[-1,1]$ onto the exterior of the unit circle, elliptic contours $\mathcal{E}_{\rho}$ with foci at points $\pm 1$ and semi-axes given by $\displaystyle \frac{1}{2}\,(\rho + \rho^{-1})\,$ and $\displaystyle \frac{1}{2}\,(\rho - \rho^{-1})\,,$ with $\rho >1$, are often considered to get suitable estimations of the error of quadrature.  Namely, these elliptic level curves are given by the expression
\begin{equation}\label{elliptic}
\mathcal{E}_{\rho} = \{z\in \mathbb{C}: |\phi(z)| = |z+\sqrt{z^2-1}| = \rho\}\,,
\end{equation}
where $\rho>1$ and the branch of $\sqrt{z^2-1}$ is taken so that $|\phi(z)|>1$ for $|z|>1$.

The outline of the paper is as follows. In Section 2, explicit expressions for the kernel \eqref{kernel} for the four induced Chebyshev weights $\displaystyle d\widehat{\sigma}_n^{[i]}\,,\,i=1,2,3,4\,,$ are provided, in such a way that they are useful to get appropriate bounds for the corresponding errors of quadrature in Section 3, which represent the main contribution of the paper. The accuracy of these bounds is checked in the fourth section by means of some illustrative numerical examples. Finally, Section 5 is devoted to the proof of the main results.

The problem of estimating the quadrature error for Gauss--type rules has been thoroughly studied in the literature; to only cite a few, see the references \cite{gauvar}, \cite{MS2003}--\cite{MS2005BIT}, \cite{Notaris2006}--\cite{Not3}, and \cite{NM2016PejSpa}--\cite{AMC2012SPP}.

\section{Explicit expressions of the kernel for the four Chebyshev weights}

The aim of this section is getting explicit expressions of the kernel $K_m$ in \eqref{kernel} corresponding to the induced measures $\displaystyle d\widehat{\sigma}_n = \pi_n^2\,d\sigma\,,$ in the particular case of the four Chebyshev weights, namely
\begin{equation*}\label{inducedweights}
\begin{array}{ll}
d\sigma^{[1]}(t) = \displaystyle \frac{dt}{\sqrt{1-t^2}}\,,\;& d\sigma^{[2]}(t) =\displaystyle  \sqrt{1-t^2}\,dt\,,\\
d\sigma^{[3]}(t) =\displaystyle   \sqrt{\frac{1-t}{1+t}}\,dt\,,\;& d\sigma^{[4]}(t) =\displaystyle  \sqrt{\frac{1+t}{1-t}}\,dt\,.
\end{array}
%\begin{split}
%d\sigma^{[1]}(t) = & \frac{dt}{\sqrt{1-t^2}}\,,\;d\sigma^{[2]}(t) = \sqrt{1-t^2}\,dt\,,\\ d\sigma^{[3]}(t) = & %\sqrt{\frac{1-t}{1+t}}\,dt\,,\;d\sigma^{[4]}(t) = \sqrt{\frac{1+t}{1-t}}\,dt\,.
%\end{split}
\end{equation*}
These explicit expressions will be used to compute different bounds for the error of quadrature, which are our main results and will be displayed in the next section. To do it, we make use of the results in \cite[\S 3]{gauli} about the explicit representations of the corresponding induced orthogonal polynomials $\displaystyle \{\widehat{\pi}_{m,n}^{[i]}\}$. Due to these results, the expression for the case of the Chebyshev weight of the first kind will be exhibited separately for $n>1$ and $n=1$; in order to distinguish both cases, hereafter $K_m^{[1]}$ will denote the corresponding kernel for $n>1$, while the other one will be referred to as $K_m^{[I]}$. While in the case where $i=1$ and $n=1$, we compute the corresponding kernel for arbitrary $m\in\mathbb N$, for the rest of the cases only the ``diagonal'' case, i.e. $m=n$, is considered for the sake of simplicity; indeed, it is also possible to deal with the general ``non--diagonal'' case $m\neq n$, but the computations are more involved and we prefer to leave it for a forthcoming paper. Therefore, specifically, the problems handled throughout this paper are the following:

\begin{itemize}

\item Computation of integrals of the form $\displaystyle I_{\omega}^{[I]}(f) = \,\int_{-1}^1\,f(t)\,\frac{t^2 dt}{\sqrt{1-t^2}}$ by means of Gauss quadrature formulae
$\displaystyle \sum_{j=1}^m\,A_{m,j}\,f(t_{m,j})\,,\,m=1,2,\ldots\,,$ where $t_{m,j}$ are the zeros of the corresponding induced orthogonal polynomials $\widehat{\pi}_{m,1}^{[1]}\,.$

\item The same for integrals of the form $\displaystyle I_{\omega}^{[i]}(f) = \,\int_{-1}^1\,f(t)\,\pi_{n,n} (t)^2\,d\sigma^{[i]}(t) dt\,,$ making use of Gauss rules with $n$ nodes taken as the zeros of the induced orthogonal polynomials $\widehat{\pi}_{n,n}^{[i]}\,,\,i=1,2,3,4\,.$

\end{itemize}

Finally, bearing in mind the well-known Joukowsky transform, the notation $\displaystyle z = \frac{1}{2}\,\left(\zeta + \frac{1}{\zeta}\right)$, $|\zeta|>1$, will be used.

Next, our conclusions are gathered in the following
\begin{lemma}\label{lem:explicit}
The explicit expression of the kernel $K_m^{[i]}$ for the four Chebyshev weights is given as follows.
\begin{itemize}

\item[\textbf{(1)}] For $i=1$ and $n=1$,
\begin{equation}\label{JezgroAp}
K_m^{[I]}(z)=\displaystyle\pi \frac{\left(\zeta^2+1\right)^2\left(1+(-1)^{m/2}\zeta^m\right)}{\zeta^{2}\left(\zeta-\zeta^{-1}\right)\left(\zeta^{2m+2}+1\right)}
\end{equation}
if $m$ is even, while
\begin{equation}\label{JezgroAn}
K_m^{[I]}(z)=\displaystyle\pi\frac{(m+2)\zeta^2+m}{\zeta^{m+2}\left(\zeta-\zeta^{-1}\right)\left(\sum_{j=0}^{(m-1)/2} (-1)^j{(m-2j)}\zeta^{m-2j}+\sum_{j=0}^{(m-1)/2} (-1)^j{(m-2j)}\zeta^{2j-m}\right)}
\end{equation}
if $m$ is odd.

\item[\textbf{(2)}] For $i=1$ and $n>1$,
\begin{equation}\label{Jezgro2}
K^{[1]}_n(z)=\displaystyle\frac{\pi\left(3\zeta^{2n}+1\right)}{2^{2n-2}\zeta^{3n}
\left(\zeta-\zeta^{-1}\right)\left(\zeta^n+\zeta^{-n}\right)}.
\end{equation}

\item[\textbf{(3)}] For $i=2$,
\begin{equation}\label{Jezgro3}
K^{[2]}_n(z)=\displaystyle\frac{\pi\left(2\xi^{2n+2}-\xi^{2n}-1\right)}{2^{2n}\xi^{3n+2}
\left(\xi-\xi^{-1}\right)\left(\xi^n+\xi^{-n}\right) }.
\end{equation}

\item[\textbf{(4)}] For $i=3$,
\begin{equation}\label{Jezgro4}
K^{[3]}_n(z)=\displaystyle\frac{\pi\left(2\zeta^{2n+1}+\zeta^{2n}+1\right)}{2^{2n}\zeta^{3n+1}
\left(\zeta-\zeta^{-1}\right)\left(\zeta^n+\zeta^{-n}\right) }.
\end{equation}

\end{itemize}
\end{lemma}

\textbf{Proof}
For the proof of this lemma, the explicit expressions for the corresponding orthogonal polynomials found in \cite[\S 3]{gauli} will be very useful. On the sequel, the orthogonal polynomials will be always monic.

\begin{itemize}

\item[\textbf{(1)}] First, in the case where $i=1$ and $n=1$, we have that $\displaystyle d\sigma_1^{[I]}=\frac{t^2\,dt}{\sqrt{1-t^2}}\,,$ and the kernel is given by $K_m^{[I]}=K^{[I]}_m(z)=\displaystyle
\frac{\varrho^{[I]}_{m}(z)}{\pi^{[I]}_{m,1}(z)},\ z\notin[-1,1]$,
where $\varrho^{[I]}_m(z)=\displaystyle\int_{-1}^1\frac{\pi^{[I]}_{m,1}(t)}{z-t}\frac{t^2\,dt}{\sqrt{1-t^2}}\,.$ Now, making the change $t=\cos{\theta}$ and using \cite[3.6]{gauli}, we have
\[
\begin{array}{rl}
\varrho^{[I]}_m(z)=&\displaystyle\int_0^{\pi}\frac{\left(\sum_{j=0}^{m/2} (-1)^j\frac{4^{-j}}{2^{m-2j-1}}\cos{(m-2j)\theta}\right)\cos^2{\theta}}{z-\cos{\theta}}\,d\theta\\[0.1in]
=&\displaystyle\frac{1}{2^{m-1}}\sum_{j=0}^{m/2} (-1)^j\int_0^{\pi}\frac{\cos{(m-2j)\theta}\,\cos^2{\theta}}{z-\cos{\theta}}\,d\theta
\end{array}
\]
if $m$ is even, and
\[
\begin{array}{rl}
\varrho^{[I]}_m(z)=&\displaystyle\int_0^{\pi}\frac{\left(\sum_{j=0}^{(m-1)/2} (-1)^j\frac{4^{-j}}{2^{m-2j-1}}\frac{m-2j}{m}\cos{(m-2j)\theta}\right)\cos^2{\theta}}{z-\cos{\theta}}\,d\theta\\[0.1in]
=&\displaystyle\frac{1}{2^{m-1}}\sum_{j=0}^{(m-1)/2} (-1)^j\frac{m-2j}{m}\int_0^{\pi}\frac{\cos{(m-2j)\theta}\,\cos^2{\theta}}{z-\cos{\theta}}\,d\theta
\end{array}
\]
if $m$ is odd. Further, using
\[
\cos{k\theta}\cos^2{\theta}=\frac{1}{2}\cos{k\theta}\left(1+\cos{2\theta}\right)=\dfrac{1}{4}\left(2\cos{k\theta}+\cos{(k-2)\theta}+\cos{(k+2)\theta} \right),
\]
and
\begin{equation}\label{int1}
\displaystyle\int_0^\pi
\frac{\cos{n\theta}}{z-\cos\theta}\,d\theta=\frac{2\pi}{\zeta^{n}(\zeta-\zeta^{-1})}, \quad n\in\mathbb{N}_0, \quad z=\frac{1}{2}\left(\zeta+\frac{1}{\zeta}\right)
\end{equation}
we get
\[
\begin{array}{rl}
\varrho^{[I]}_m(z)=&\displaystyle\frac{1}{2^{m-1}}\sum_{j=0}^{m/2} (-1)^j\dfrac{1}{4}\frac{2\pi}{\zeta-\zeta^{-1}}\left(2\zeta^{2j-m}+\zeta^{-|2j+2-m|}+\zeta^{2j-m-2}\right)\\
=&\displaystyle\pi\frac{\sum_{j=0}^{m/2}(-1)^j\left(2\zeta^{2j-m}+\zeta^{-|m-2j-2|}+\zeta^{2j-m-2}\right)}{{2^m}\left(\zeta-\zeta^{-1}\right)} \\ = &\displaystyle\pi\frac{\zeta^{-m-2}+\zeta^{-m}+(-1)^{m/2}\left(\zeta^{-2}+1\right)}{{2^m}\left(\zeta-\zeta^{-1}\right)}=\pi \frac{\left(\zeta^2+1\right)\left(1+(-1)^{m/2}\zeta^m\right)}{2^m\zeta^{m+2}}
\end{array}
\]
if $m$ is even
and
\[
\begin{array}{rl}
\varrho^{[I]}_m(z)=&\displaystyle\frac{1}{2^{m-1}}\sum_{j=0}^{(m-1)/2} (-1)^j(m-2j)\dfrac{1}{4}\frac{2\pi}{\zeta-\zeta^{-1}}\left(2\zeta^{2j-m}+\zeta^{-|m-2j-2|}+\zeta^{2j-m-2}\right)\\
=&\displaystyle\pi\frac{\sum_{j=0}^{(m-1)/2}(-1)^j(m-2j)\left(2\zeta^{2j-m}+\zeta^{-|m-2j-2|}+\zeta^{2j-m-2}\right)}{m{2^m}\left(\zeta-\zeta^{-1}\right)} \\ =&\displaystyle\pi\frac{m\zeta^{-m-2}+(m+2)\zeta^{-m}}{ m{2^m}\left(\zeta-\zeta^{-1}\right)}
\end{array}
\]
if $m$ is odd. Now, using the substitutions
\begin{equation}\label{Ceb1}
z=\frac{1}{2}\left(\zeta+\frac{1}{\zeta}\right), \quad T_j\left(\frac{1}{2}\left(\zeta+\frac{1}{\zeta}\right)\right)=\frac{1}{2}\left(\zeta^j+\frac{1}{\zeta^j}\right),
\end{equation}
and \cite[3.6]{gauli} again, we get
\[
\begin{array}{rl}
\pi^{[I]}_{m,1}(z)=&\displaystyle\frac{1}{2^{m-1}}\sum_{j=0}^{m/2} (-1)^j\frac{1}{2}\left(\zeta^{m-2j}+\frac{1}{\zeta^{m-2j}}\right)\\
=&\displaystyle\frac{\sum_{j=0}^{m/2} (-1)^j\zeta^{m-2j}+\sum_{j=0}^{m/2} (-1)^j\zeta^{2j-m}}{2^m}=\displaystyle\dfrac{\zeta^{2m+2}+1}{2^m\zeta^m\left(\zeta^2+1\right)}
\end{array}
\]
if $m$ is even and \eqref{JezgroAp} is established. In the same way,
\[
\pi^{[I]}_{m,1}(z)=\displaystyle\frac{1}{2^{m-1}}\sum_{j=0}^{(m-1)/2} (-1)^j\frac{m-2j}{m}\frac{1}{2}\left(\zeta^{m-2j}+\frac{1}{\zeta^{m-2j}}\right)
\]
when $m$ is odd and \eqref{JezgroAn} holds. The last sum also can be explicitly calculated, but it is not simple expresion as in the case of even number $m$.

\item[\textbf{(2)}] In this case (see \cite[3.4]{gauli}), $\displaystyle d\sigma_n^{[1]}=\frac{\mathring{T}_n^2(t)\,dt}{\sqrt{1-t^2}}$ and the kernel is given by $K^{[1]}_n(z)=\displaystyle
\frac{\varrho^{[1]}_{n}(z)}{\pi^{[1]}_{n,n}(z)},\ z\notin[-1,1]$,
where $\pi^{[1]}_{n,n}(t)=\mathring{T}_n(t)$, $\varrho^{[1]}_n(z)=\displaystyle\int_{-1}^1\frac{\mathring{T}_n(t)}{z-t}\frac{\mathring{T}_n^2(t)\,dt}{\sqrt{1-t^2}}$ and $\mathring{T}_n(t)=\frac{1}{2^{n-1}}T_n(t)\,$ (that is, the monic Chebyshev polynomial of the first kind), which means
\[
\begin{array}{rl}
\varrho^{[1]}_n(z)=&\displaystyle\frac{1}{2^{3n-3}}\displaystyle\int_0^{\pi}\frac{\cos^3{n\theta}}{z-\cos{\theta}}\,d\theta
\end{array}
\]
after we put $t=\cos{\theta}$. Further, using
\[
\cos^3{n\theta}=\dfrac{3\cos{n\theta}+\cos{3n\theta}}{4},
\]
and (\ref{int1}),
we get
\begin{equation}\label{ro2}
\begin{array}{rl}
\varrho^{[1]}_n(z)=&\displaystyle\frac{1}{2^{3n-1}}\left(\frac{6\pi}{\zeta^{n}(\zeta-\zeta^{-1})}+\frac{2\pi}{\zeta^{3n}(\zeta-\zeta^{-1})} \right)
\end{array}
\end{equation}
and using again (\ref{Ceb1}), the proof of \eqref{Jezgro2} is fulfilled.

\item[\textbf{(3)}] Now, $\displaystyle d\sigma_n^{[2]}={\mathring{U}_n^2(t){\sqrt{1-t^2}}\,dt}\,$ (see \cite[Theorem 3.4]{gauli}), where $\mathring{U}_n$ denotes the monic Chebyshev polynomial of the second kind, and the kernel is given by $K^{[2]}_n(z)=\displaystyle
\frac{\varrho^{[2]}_{n}(z)}{\pi^{[2]}_n(z)},\ z\notin[-1,1]$,
where $\pi^{[2]}_{n,n}(t)=\mathring{T}_n(t)$, $\varrho^{[2]}_n(z)=\displaystyle\int_{-1}^1\frac{\mathring{T}_n(t)}{z-t}{\mathring{U}_n^2(t)}{\sqrt{1-t^2}}\,dt$ and $\mathring{U}_n(t)=\frac{1}{2^{n}}U_n(t)$, which means
\[
\begin{array}{rl}
\varrho^{[2]}_n(z)=&\displaystyle\frac{1}{2^{3n-1}}\displaystyle\int_0^{\pi}\frac{\cos{n\theta}\sin^2{(n+1)\theta}}{z-\cos{\theta}}\,d\theta\,.
\end{array}
\]
Further, using the identity
\begin{eqnarray*}
\cos{n\theta}\sin^2{(n+1)\theta}&=&\frac{1}{2}\cos{n\theta}\left(1-\cos{2(n+1)\theta}\right)\\
&=&\dfrac{1}{4}\left(2\cos{n\theta}-\cos{(n+2)\theta}-\cos{(3n+2)\theta} \right),
\end{eqnarray*}
and (\ref{int1}),
we get
\begin{equation}\label{ro3}
\begin{array}{rl}
\varrho^{[2]}_n(z)=&\displaystyle\frac{1}{2^{3n-1}}\dfrac{1}{4}\frac{2\pi}{\zeta-\zeta^{-1}}\left(\frac{2}{\zeta^{n}}-\frac{1}{\zeta^{n+2}}-\frac{1}{\zeta^{3n+2}}\right)
\end{array}
\end{equation}
and after using (\ref{Ceb1}) again, we have that \eqref{Jezgro3} holds.

\item[\textbf{(4)}] In this case, $\displaystyle d\sigma_n^{[3]}={\mathring{V}_n^2(t)}{\sqrt{\frac{1+t}{1-t}}}\,dt\,$ (see \cite[Theorem 3.6]{gauli}), where $\mathring{V}_n$ denotes the monic Chebyshev polynomial of the third kind, and the kernel is given by $K^{[3]}_n(z)=\displaystyle
\frac{\varrho^{[3]}_{n}(z)}{\pi^{[3]}_n(z)},\ z\notin[-1,1]$,
where $\pi^{[3]}_{m,n}(t)=\mathring{T}_n(t)$, $\varrho^{[3]}_n(z)=\displaystyle\int_{-1}^1\frac{\mathring{T}_n(t)}{z-t}{\mathring{V}_n^2(t)}{\sqrt{\frac{1+t}{1-t}}}\,dt$ and $\mathring{V}_n(t)=\frac{1}{2^{n}}V_n(t)$, which implies
\[
\begin{array}{rl}
\varrho^{[3]}_n(z)=&\displaystyle\frac{1}{2^{3n-1}}\displaystyle\int_0^{\pi}\frac{\cos{n\theta}\cos^2{(n+1/2)\theta}}{z-\cos{\theta}}\,d\theta.
\end{array}
\]
Then, proceeding analogously as above and using the identity
\begin{eqnarray*}
\cos{n\theta}\cos^2{(n+1/2)\theta}&=&\frac{1}{2}\cos{n\theta}\left(1+\cos{(2n+1)\theta}\right)\\
&=&\dfrac{1}{4}\left(2\cos{n\theta}+\cos{(n+1)\theta}+\cos{(3n+1)\theta} \right),
\end{eqnarray*}
and (\ref{int1}),
we get
\begin{equation}\label{ro4}
\begin{array}{rl}
\varrho^{[3]}_n(z)=&\displaystyle\frac{1}{2^{3n-1}}\dfrac{1}{4}\frac{2\pi}{\zeta-\zeta^{-1}}\left(\frac{2}{\zeta^{n}}+\frac{1}{\zeta^{n+1}}+\frac{1}{\zeta^{3n+1}}\right)
\end{array}
\end{equation}
and after making use of (\ref{Ceb1}), \eqref{Jezgro4} is easily proven.

\end{itemize}

\begin{remark}\label{rem:34}
We have not shown the results for $i=4$ in Lemma \ref{lem:explicit}, since the orthogonal polynomials $\displaystyle \widehat{\pi}_{k,n}^{[4]}$ and $\displaystyle \widehat{\pi}_{k,n}^{[3]}$ are easily related to each other. Indeed (see e.g. \cite[(3.15)]{gauli},
$$\widehat{\pi}_{k,n}^{[4]}(t) = (-1)^k\,\widehat{\pi}_{k,n}^{[3]}(-t)\,$$
and thus, it is enough to consider the case for $i=3$.
\end{remark}

\section{Main results}

Next, using the results in previous Lemma \ref{lem:explicit}, different bounds of the error of quadrature for the four induced quadrature weights are derived. They are presented in the following subsections.

\subsection{$L^{\infty}$--bounds}

Hereafter, for a function $g$ and a compact subset $E$ of the complex plane, the $L^{\infty}$--norm of $g$ on $E$ will be denoted by
$$\|g\|_E = \max_{z\in E}\,|g(z)|\,.$$ Now, from \eqref{intrepres}, taking $\Gamma = \mathcal{E}_{\rho}\,,$ we easily get that if $f$ is analytic on $\mathcal{E}_{\rho}$ and its interior, for certain $\rho >1$, then,
\begin{equation}\label{errorbound}
|R_m(f)| \leq \,\frac{l(\mathcal{E}_{\rho})}{2\pi}\,\|K_m\|_{\mathcal{E}_{\rho}}\,\|f(z)\|_{\mathcal{E}_{\rho}}\,.
\end{equation}
On the sequel, if we denote by $D_{\rho}$ the closed interior of $\mathcal{E}_{\rho}$, define
$$\rho_f = \sup \{\rho >1\,: f\;\;\text{is analytic on}\;\;D_{\rho}\}\,.$$
Next, we are concerned with the maximum modulus of the kernel corresponding to each of the cases considered in Lemma \ref{lem:explicit} on the level curves $\mathcal{E}_{\rho}$ defined in \eqref{elliptic}. The results are shown in the following theorem.
\begin{theorem}\label{thm:maxmod}
The maximum modulus of the kernel for the four Chebyshev weights is given as follows.

\begin{itemize}

\item[\textbf{(1)}] For each $m\in\mathbb N$ there exists $\rho^*=\rho^*(m)>1$ such that for each $\rho >\rho^*$ the modulus of the kernels (\ref{JezgroAp}) and (\ref{JezgroAn}) attain their maximum value on the real axis, i.e. \\[0.1in]
\[
\|K^{[I]}_m\|_{\mathcal{E}_{\rho}} = \left|K_{m}^{[I]}\left(-{\textstyle\frac{1}{2}}(\rho+\rho^{-1})\right)\right|=
\left|K^{[I]}_{m}\left({\textstyle\frac{1}{2}}(\rho+\rho^{-1})\right)\right|.
\]

\item[\textbf{(2)}] Let $i=1$ and $n>1$. Then, for each $n\in\mathbb N$ there exists $\rho^*=\rho^*(n)>1$ such that for each $\rho >\rho^*$ the modulus of the kernel (\ref{Jezgro2}) attains its maximum value on the real axis, i.e. \\[0.1in]
\[
\|K^{[1]}_n\|_{\mathcal{E}_{\rho}} = \left|K^{[1]}_{n}\left(-{\textstyle\frac{1}{2}}(\rho+\rho^{-1})\right)\right|=
\left|K^{[1]}_{n}\left({\textstyle\frac{1}{2}}(\rho+\rho^{-1})\right)\right|.
\]

\item[\textbf{(3)}] Let $i=2$. For each $n\in\mathbb N/\{1\}$ there exists $\rho^*=\rho^*(n)>1$ such that for each $\rho >\rho^*$ the modulus of the kernel (\ref{Jezgro3}) attains its maximum value on the real axis,
while for $n=1$ there exists $\rho^*=\rho^*(n)>1$ such that for each $\rho >\rho^*$ the modulus of the kernel (\ref{Jezgro3}) attains its maximum value on the imaginary axis, i.e. \\[0.1in]
\[
\|K^{[2]}_n\| = \left|K^{[2]}_{n}\left(-{\textstyle\frac{1}{2}}(\rho+\rho^{-1})\right)\right|=
\left|K^{[2]}_{n}\left({\textstyle\frac{1}{2}}(\rho+\rho^{-1})\right)\right|,
\]
for $n>1$ and
\[
\|K^{[2]}_1\|_{\mathcal{E}_{\rho}} = \left|K^{[2]}_{1}\left(\frac{i}{2}(\rho-\rho^{-1})\right)\right|.
\]

\item[\textbf{(4)}] Let $i=3$. For each $n\in\mathbb N$ there exists $\rho^*=\rho^*(n)>1$ such that for each $\rho >\rho^*$ the modulus of the kernel (\ref{Jezgro4}) attains its maximum value on the real axis, i.e. \\[0.1in]
\[
\|K^{[3]}_n\|_{\mathcal{E}_{\rho}} = \left|K^{[3]}_{n}\left(-{\textstyle\frac{1}{2}}(\rho+\rho^{-1})\right)\right|=
\left|K^{[3]}_{n}\left({\textstyle\frac{1}{2}}(\rho+\rho^{-1})\right)\right|.
\]

\end{itemize}

\end{theorem}

We empirically found that the values $\rho^*$ in the previous theorem are relatively closed to $1$ in all the cases, which gives us an opportunity to successfully apply the estimate based on the maximum modulus of the kernel. Anyway, hereafter we assume that for the integrand $f$, it holds $\rho^* \leq \rho_f\,.$ Indeed, since the length of the ellipse can be estimated by
\begin{equation*}
 l(\mathcal{E}_{\rho})\leqslant 2\pi
a_1\left(1-\frac{1}{4}a_1^{-2}-\frac{3}{64}a_1^{-4}-\frac{5}{256}a_1^{-6}\right)\,,
\end{equation*}
where $a_1=\frac{1}{2}(\rho+\rho^{-1})$, by \eqref{errorbound} and the results in Theorem \ref{thm:maxmod}, the $L^{\infty}$--bounds for the error of quadrature may be written as follows.
\begin{equation}\label{prvaIparno}
r^{[I]}_1(f)=\inf_{\rho^*<\rho<\rho_f} \displaystyle\pi a_1 \frac{\left(\rho^2+1\right)^2\left(1+(-1)^{m/2}\rho^m\right)\left(1-\frac{1}{4}a_1^{-2}-\frac{3}{64}a_1^{-4}-\frac{5}{256}a_1^{-6}\right)\,\|f\|_{\mathcal{E}_{\rho}}}{\rho^{2}\left(\rho-\rho^{-1}\right)\left(\rho^{2m+2}+1\right)}
\end{equation}
if $m$ is even, and
\begin{equation}\label{prvaIneparno}
r^{[I]}_1(f)=\inf_{\rho^*<\rho<\rho_f} \frac{\pi a_1\left((m+2)\rho^2+m\right)\left(1-\frac{1}{4}a_1^{-2}-\frac{3}{64}a_1^{-4}-\frac{5}{256}a_1^{-6}\right)\,\|f\|_{\mathcal{E}_{\rho}}
	}{\rho^{m+2}\left(\rho-\rho^{-1}\right)\left(\sum_{j=0}^{(m-1)/2} (-1)^j{(m-2j)}\rho^{m-2j}+\sum_{j=0}^{(m-1)/2} (-1)^j{(m-2j)}\rho^{2j-m}\right)}
\end{equation}
if $m$ is odd. In the same way,
\begin{equation}\label{prva1}
r^{(1)}_1(f)=\inf_{\rho^*<\rho<\rho_f} \displaystyle\frac{\pi a_1\left(3\rho^{2n}+1\right)\left(1-\frac{1}{4}a_1^{-2}-\frac{3}{64}a_1^{-4}-\frac{5}{256}a_1^{-6}\right)\,\|f\|_{\mathcal{E}_{\rho}}
	}{2^{2n-2}\rho^{3n}\left(\rho-\rho^{-1}\right)\left(\rho^n+\rho^{-n}\right)},
\end{equation}
\begin{equation}\label{prva2}
r^{(2)}_1(f)=\inf_{\rho^*<\rho<\rho_f} \displaystyle\frac{\pi a_1\left(2\rho^{2n+2}-\rho^{2n}-1\right)\left(1-\frac{1}{4}a_1^{-2}-\frac{3}{64}a_1^{-4}-\frac{5}{256}a_1^{-6}\right)\,\|f\|_{\mathcal{E}_{\rho}}
	}{2^{2n}\rho^{3n+2}\left(\rho-\rho^{-1}\right)\left(\rho^n+\rho^{-n}\right)}, \quad n>1,
\end{equation}
\begin{equation}\label{treca1}
r^{(3)}_1(f)=\inf_{\rho^*<\rho<\rho_f}\displaystyle\frac{\pi a_1\left(2\rho^{2n+1}+\rho^{2n}+1\right)\left(1-\frac{1}{4}a_1^{-2}-\frac{3}{64}a_1^{-4}-\frac{5}{256}a_1^{-6}\right)\,\|f\|_{\mathcal{E}_{\rho}}
	}{2^{2n}\rho^{3n+1}\left(\rho-\rho^{-1}\right)\left(\rho^n+\rho^{-n}\right)},
\end{equation}
where $\rho^*$ is the value defined above, which has been obtained empirically.

\subsection{Error bounds based on the expansion of the remainder}

If $f$ is an analytic function in a neighborhood of the real interval $[-1,1]$, there exists $\rho > 1$ such that $f$ is analytic in the interior of
$\mathcal{E}_{\rho}$ and, so, it admits the expansion
\begin{equation}\label{cebred}
f(z)=\sum_{k=0}^{\infty}{'}\alpha_k T_k(z),
\end{equation}
where $\alpha_k$ are given by
\[
\alpha_k=\dfrac{1}{\pi}\int_{-1}^1(1-t^2)^{-1/2}f(t)T_k(t)dt.
\]
The series (\ref{cebred}) converges for each $z$ in the interior
of $\mathcal{E}_{\rho}$. The ``prim'' symbol in the corresponding sum denotes
that the first term is taken with the factor ${1}/{2}$. In general, the Chebyshev-Fourier coefficients $\alpha_k$ in
(\ref{cebred}) are unknown.
However, Elliott \cite{Eliot} described a number of ways of estimating or bounding them. In particular, for our purpose it is useful the upper bound
\begin{equation}\label{cebocena}
|\alpha_k|\leq \dfrac{2}{\rho^k}\,\|f\|_{\mathcal{E}_{\rho}}.
\end{equation}
Next, the bounds for the error of quadrature obtained by expanding the remainder terms \eqref{errorbound} are listed. The case corresponding to $i=1$ and $n=1$ is omitted since the computations are too complicated.
\begin{theorem}\label{thm:errors2}

The following upper bounds for the error of quadrature hold.

\begin{itemize}

\item[\textbf{(1)}] For $i=1$ and $n>1$,
\begin{equation}\label{Druga1}
r^{(1)}_2(f)=\inf_{1<\rho<\rho_f} \dfrac{\pi}{2^{2n-2}}\dfrac{1}{\rho^{2n}-1}\,\|f\|_{\mathcal{E}_{\rho}}.
\end{equation}

\item[\textbf{(2)}] For $i=2$,
\begin{equation}\label{Druga2}
r^{(2)}_2(f)=\inf_{1<\rho<\rho_f}\dfrac{\pi}{2^{2n}}\left(\dfrac{1}{\rho^{2n}-1}+\dfrac{1}{2\rho^{2n+2}}\right)\,\|f\|_{\mathcal{E}_{\rho}}.
\end{equation}

\item[\textbf{(3)}] For $i=3$,
\begin{equation}\label{Druga3}
	r^{(3)}_2(f)=\inf_{1<\rho<\rho_f}\dfrac{\pi}{2^{2n}}\left(\dfrac{1}{\rho^{2n}-1}\right)\,\|f\|_{\mathcal{E}_{\rho}}.
	\end{equation}

\end{itemize}

\end{theorem}

\subsection{$L^{1}$--bounds}

In this subsection, our goal consists in bounding the integrals

\begin{equation}\label{L1int}
L_m^{[i]}(\mathcal{E}_{\rho})=\dfrac{1}{2\pi}\oint_{\mathcal{E}_{\rho}}\left|K_m^{[i]}(z)\right||dz|\,,
\end{equation}
in the case of the four Chebyshev weights considered. In the case of $i=1$ and $n=1$, the computations with the kernel $K_m^{[I]}(z)$ are again quite involved and so, we prefer to omit the results. Otherwise, we have the following upper bounds for \eqref{L1int}.
\begin{theorem}\label{thm:errors3}

We have the following $L^1$--error bounds.

\begin{itemize}

\item[\textbf{(1)}] For $i=1$ and $n>1$,
\begin{equation}\label{Treca1}
r^{(1)}_3(f)=\inf_{1<\rho<\rho_f} \dfrac{\pi}{\rho^n\cdot 2^{2n-1}}\sqrt{\dfrac{7\rho^{-2n}+9\rho^{2n}}{\rho^{4n}-1}}\,\|f\|_{\mathcal{E}_{\rho}}.
\end{equation}

\item[\textbf{(2)}] For $i=2$,
\begin{equation}\label{Treca2}
r^{(2)}_3(f)=\inf_{1<\rho<\rho_f} \dfrac{\pi}{\rho^n\cdot 2^{2n+1}}\sqrt{\dfrac{\rho^{2n-4}+4\rho^{2n}+3\rho^{-2n-4}}{\rho^{4n}-1}}\,\|f\|_{\mathcal{E}_{\rho}}.
\end{equation}

\item[\textbf{(3)}] For $i=3$,
\begin{equation}\label{Treca3}
r^{(3)}_3(f)=\inf_{1<\rho<\rho_f}\dfrac{\pi}{\rho^n\cdot 2^{2n+1}}\sqrt{\dfrac{\rho^{2n-2}+4\rho^{2n}+3\rho^{-2n-2}}{\rho^{4n}-1}}\,\|f\|_{\mathcal{E}_{\rho}}.
\end{equation}

\end{itemize}

\end{theorem}

\section{Numerical experiments}

Throughout this section, we test the previously given error bounds, i.e. \eqref{prvaIparno}--\eqref{treca1}, \eqref{Druga1}--\eqref{Druga3} and \eqref{Treca1}--\eqref{Treca3}, by means of the two characteristic examples: $f_0(z)=e^{\omega z^2}, \ \omega>0$ and $f_1(z)=e^{\cos (\omega z)},\ \omega>0$. It is clear that both $f_0$ and $f_1$ are entire functions, so that $\rho_f = +\infty\,.$
We know that
\[
\|f_0\|_{\mathcal{E}_{\rho}}=e^{\omega a_1^2}, \quad \|f_1\|_{\mathcal{E}_{\rho}}=e^{\cosh{\omega b_1} },
\]
where $a_1=\frac{\rho+\rho^{-1}}{2}$ and $b_1=\frac{\rho-\rho^{-1}}{2}$.
The corresponding results for both test-functions concerned with each of the considered weights are presented in the Tables \ref{figura0a} -- \ref{figura1g}. In the case $i=1$ and $n=1$, corresponding to Tables \ref{figura0a} and \ref{figura1a}, only odd values of $m$ are considered since the estimate (\ref{prvaIneparno}) gives much better results than (\ref{prvaIparno}); this seems natural as we look the leading coeeficients of $\rho$ in the numerator and the denominator in those expressions.

It is noteworthy that in general the estimates of the error are quite sharp, as well as the accuracy of the respective quadrature rules.

\begin{table}[!htp]
	\begin{center}{\footnotesize
			\begin{tabular}{llll}\hline
				$m,\omega$ & $r^{[I]}_1(f_0)$ & ${\rm Error}^{[I]}$ & $I_{\omega}^{[I]} (f_0)$ \\
				\hline
				$5, 0.5$ & $4.351(-6)$ &  $ 7.398(-7)$ &
				$2.3026...(+0)$ \\
				$25, 0.5$ & $4.747(-47)$ &  $ 3.753(-48)$ &
				$2.3026...(+0)$ \\
				$35, 0.5$ & $ 7.646(-71)$ &  $5.122(-72)$ &
				$2.3026...(+0)$ \\
				\hline
		$5, 1$ & $   1.912(-4)$ &  $ 3.147(-5)$ &
		$3.4221...(+0)$ \\
		$25, 1$ & $2.077(-39)$ &  $ 1.634(-40)$ &
		$3.4221...(+0)$ \\
		$35, 1$ & $3.414(-60)$ &  $2.278(-61)$ &
		$3.4221...(+0)$ \\
		\hline
			$5, 5$ & $9.162(+0)$ &  $ 1.130(+0)$ &
		$1.1111...(+2)$ \\
		$25, 5$ & $ 5.439(-21)$ &  $ 4.084(-22)$ &
		$1.1111...(+2)$ \\
		$35, 5$ & $8.313(-35)$ &  $5.376(-36)$ &
		$1.1111...(+2)$ \\
		\hline
	\end{tabular}}
		\vskip 0.2in
		\caption{The values of the bounds $r_1^{[I]}(f_0)$, compared with the actual value of the error ${\rm Error}^{[I]}$ and the value of the integral $I_{\omega}^{[I]} (f_0)$ for some values of $m$, $\omega$ in the case of $d\sigma_n^{[I]}$}\label{figura0a}
	\end{center}
\end{table}

\begin{table}[!htp]
	\begin{center}{\footnotesize
			\begin{tabular}{llll}\hline
				$m,\omega$ & $r^{[I]}_1(f_1)$ & ${\rm Error}^{[I]}$ & $I_{\omega}^{[I]} (f_1)$ \\
				\hline
				$5, 0.5$ & $6.901(-8)$ &  $ 8.926(-9)$ &
				$3.8958...(+0)$ \\
				$25, 0.5$ & $9.011(-51)$ &  $ 1.202(-52)$ &
				$3.8958...(+0)$ \\
				$35, 0.5$ & $8.710(-74)$ &  $3.733(-75)$ &
				$3.8958...(+0)$ \\
				\hline
				$5, 1$ & $ 5.363(-5)$ &  $ 6.051(-6)$ &
				$3.0296...(+0)$ \\
				$25, 1$ & $4.814(-36)$ &  $ 2.299(-37)$ &
				$3.0296...(+0)$ \\
				$35, 1$ & $4.091(-53)$ &  $1.632(-54)$ &
				$3.0296...(+0)$ \\
				\hline
				$5, 5$ & $5.315(+0)$ &  $ 6.236(-2)$ &
				$1.3763...(+2)$ \\
				$25, 5$ & $ 4.584(-8)$ &  $ 3.481(-10)$ &
				$1.3763...(+2)$ \\
				$35, 5$ & $ 4.274(-13)$ &  $3.242(-15)$ &
				$1.3763...(+2)$ \\
				\hline
		\end{tabular}}
		\vskip 0.2in
		\caption{The values of the bounds $r_1^{[I]}(f_1)$, compared with the actual value of the error ${\rm Error}^{[I]}$ and the value of the integral $I_{\omega}^{[I]} (f_1)$ for some
			values of $m$, $\omega$ in the case of $d\sigma_n^{[I]}$}\label{figura1a}
	\end{center}
\end{table}

\begin{table}[!htp]
	\begin{center}{\footnotesize
			\begin{tabular}{llllll}\hline
				$n,\omega$ & $r^{[1]}_1(f_0)$ & $r^{[1]}_2(f_0)$
				& $r^{[1]}_3(f_0)$ & ${\rm Error}^{[1]}$ & $I_{\omega}^{[1]} (f_0)$ \\
				\hline
				$6, 0.1$ & $1.026(-14) $ & $ 6.809(-15)$ & $3.209(-14)$  & $ 1.641(-15)$ &
				$1.6136...(-3)$ \\
				$8, 0.1$ & $8.224(-21) $ & $5.466(-21) $ &$8.199(-21)$ & $ 1.144(-21)$ &
				$1.0085...(-4)$ \\
				$10, 0.1$ & $3.980(-27) $ & $2.647(-27) $ & $3.970(-27)$ & $4.967(-28)$ &
				$6.3032...(-6)$ \\
				$15, 0.1 $ & $ 1.285(-43)$ & $8.554(-44) $ & $1.283(-43)$ & $ 1.314(-44)$&
				$6.1555...(-9)$ \\
				$20, 0.1 $ & $ 7.594(-61)$ & $5.056(-61) $ & $7.584(-61)$ & $ 6.738(-62)$&
				$6.0112...(-12)$
				\\
				\hline
				$6, 1$ & $1.690(-8) $ & $ 1.079(-8)$ & $1.618(-8)$ & $ 2.596(-9)$ &
				$2.6896...(-3)$ \\
				$8, 1$ & $1.338(-12) $ & $8.645(-13) $ &$1.297(-12)$ & $ 1.807(-13)$ &
				$1.6810...(-4)$ \\
				$10, 1$ & $6.430(-17) $ & $4.178(-17) $ & $6.267(-17)$ & $7.833(-18)$ &
				$1.0506...(-5)$ \\
				$15, 1 $ & $ 2.058(-28)$ & $1.349(-28) $ & $2.023(-28)$ & $ 2.069(-29)$ &
				$1.0260...(-8)$ \\
				$20, 1 $ & $ 1.210(-40)$ & $7.962(-41) $ & $1.194(-40)$  & $ 1.060(-41)$ &
				$1.0019...(-11)$ \\
				\hline
				$6, 5$ & $3.000(-3) $ & $ 1.600(-3)$ & $2.400(-3)$ & $ 3.702(-4)$ &
				$6.1602...(-2)$ \\
				$8, 5$ & $5.337(-6) $ & $3.013(-6) $ & $4.519(-6)$ & $ 6.153(-7)$ &
				$3.8426...(-3)$  \\
				$10, 5$ & $6.014(-9) $ & $3.519(-9) $ & $5.279(-9)$  & $6.473(-10)$ &
				$2.4015...(-4)$ \\
				$15, 5 $ & $ 5.499(-17)$ & $3.361(-17) $ & $5.041(-17)$ & $ 5.124(-18)$ &
				$2.3452...(-7)$ \\
				$20, 5 $ & $ 9.677(-26)$ & $6.042(-26) $ & $9.063(-26)$ & $ 8.020(-27)$ &
				$2.2902...(-10)$\\
				\hline
		\end{tabular}}
		\vskip 0.2in
		\caption{The values of the derived bounds $r_1^{[1]}(f_0),r_2^{[1]}(f_0),r_3^{[1]}(f_0)$, compared with the actual value of the error ${\rm Error}^{[1]}$ and the value of the integral $I_{\omega}^{[1]} (f_0)$ for some
			values of $n$, $\omega$ in the case of $d\sigma_n^{[1]}$}\label{figura0b}
	\end{center}
\end{table}

\begin{table}[!htp]
	\begin{center}{\footnotesize
			\begin{tabular}{llllll}\hline
				$n,\omega$ & $r^{[1]}_1(f_1)$ & $r^{[1]}_2(f_1)$
				& $r^{[1]}_3(f_1)$ & ${\rm Error}^{[1]}$ & $I_{\omega}^{[1]} (f_1)$ \\
				\hline
				$6, 0.2$ & $3.154(-17) $ & $ 2.098(-17)$ &$3.148(-17)$  & $ 3.840(-18)$ &
				$4.1285...(-3)$ \\
				$10, 0.2$ & $7.275(-31) $ & $4.843(-31) $ &$7.265(-31)$  & $6.599(-32)$ &
				$1.6126...(-5)$ \\
				$15, 0.2 $ & $ 1.946(-48)$ & $1.296(-48) $ &$1.944(-48)$  & $ 1.595(-49)$ &
				$1.5749...(-8)$ \\
				$20, 0.2 $ & $ 2.223(-66)$ & $1.480(-66) $ &$2.221(-66)$  & $ 1.349(-67)$ &
				$1.5379...(-11)$\\
				\hline
				$6, 1$ & $4.856(-9) $ & $ 3.095(-9)$ & $4.643(-9)$ & $ 5.596(-10)$ &
				$3.3409...(-3)$ \\
				$10, 1$ & $3.793(-17) $ & $2.444(-17) $ &$3.666(-17)$  & $3.297(-18)$ &
				$1.3050...(-5)$ \\
				$15, 1 $ & $ 8.548(-28)$ & $5.545(-28) $ &$8.317(-28)$  & $ 5.915(-29)$&
				$1.2744...(-8)$ \\
				$20, 1 $ & $ 8.371(-39)$ & $5.448(-39) $ & $8.172(-39)$ & $ 4.922(-40)$&
				$1.2446...(-11)$
				\\
				\hline
				$6, 5$ & $3.200(-3) $ & $ 9.681(-4)$ & $1.400(-4)$ & $ 1.240(-4)$ &
				$1.5738...(-3)$ \\
				$10, 5$ & $6.145(-7) $ & $2.128(-7) $ &$3.191(-7)$ & $2.264(-8)$ &
				$5.9929...(-6)$ \\
				$15, 5 $ & $ 7.284(-12)$ & $2.811(-12) $ &$4.216(-12)$ & $ 2.437(-13)$ &
				$5.8450...(-9)$ \\
				$20, 5 $ & $ 5.363(-17)$ & $2.248(-17) $ &$3.371(-17)$& $ 1.584(-18)$ &
				$5.7081...(-12)$ \\
				\hline
		\end{tabular}}
		\vskip 0.2in
		\caption{The values of the derived bounds $r_1^{[1]}(f_1),r_2^{[1]}(f_1),r_3^{[1]}(f_1)$, compared with the actual value of the error ${\rm Error}^{[1]}$ and the value of the integral $I_{\omega}^{[1]} (f_1)$ for some
			values of $n$, $\omega$ in the case of $d\sigma_n^{[1]}$}\label{figura1b}
	\end{center}
\end{table}

\begin{table}[!htp]
	\begin{center}{\footnotesize
			\begin{tabular}{llllll}\hline
				$n,\omega$ & $r^{[2]}_1(f_0)$ & $r^{[2]}_2(f_0)$
				& $r^{[2]}_3(f_0)$ & ${\rm Error}^{[2]}$ & $I_{\omega}^{[2]} (f_0)$ \\
				\hline
				$5, 0.1$ & $1.500(-12) $ &  $1.500(-12) $ & $1.496(-12)$  & $ 2.620(-13)$ &
				$1.6136...(-3)$ \\
				$10, 0.1$ & $6.625(-28) $ & $ 6.625(-28)$ & $6.617(-28)$ & $ 8.269(-29)$ &
				$1.5758...(-6)$ \\
				$15, 0.1$ & $2.140(-44) $ & $2.140(-44) $ &$2.139(-44)$  & $2.189(-45)$ &
				$1.5388...(-9)$ \\
				$20, 0.1 $ & $ 1.265(-61)$ & $ 1.265(-61)$ &$1.264(-61)$  & $ 1.222(-62)$&
				$1.5028...(-12)$ \\
				\hline
				$5, 1$ & $2.439(-7) $ & $2.435(-7) $ &$2.377(-7)$ & $ 4.073(-8)$ &
				$2.6896...(-3)$ \\
				$10, 1$ & $1.058(-17) $ & $1.058(-17) $ &$1.045(-17)$  & $1.291(-18)$ &
				$1.6266...(-6)$ \\
				$15, 1 $ & $ 3.401(-29)$ & $ 3.400(-29)$ & $3.371(-29)$ & $ 3.422(-30)$ &
				$2.5650...(-9)$ \\
				$20, 1 $ & $ 2.003(-41)$ & $ 2.003(-41)$ &  $1.991(-41)$& $ 1.756(-42)$ &
				$2.5049...(-12)$ \\
				\hline
				$5, 5$ & $8.600(-3) $ &$8.200(-3) $  &$7.400(-3)$ & $ 1.104(-3)$ &
				$6.1356...(-2)$  \\
				$10, 5$ & $9.428(-10) $ &$9.321(-10) $  &$8.814(-10)$  & $1.018(-10)$ &
				$6.0038...(-5)$ \\
				$15, 5 $ & $ 8.791(-18)$ & $ 8.745(-18)$ &$8.409(-18)$  & $ 8.206(-19)$ &
				$5.8631...(-8)$ \\
				$20, 5 $ & $ 1.562(-26)$ & $ 1.558(-26)$ &$1.511(-26)$  & $ 1.297(-27)$ &
				$2.7257...(-11)$\\
				\hline
		\end{tabular}}
		\vskip 0.2in
		\caption{The values of the derived bounds $r_1^{[2]}(f_0),r_2^{[2]}(f_0),r_3^{[2]}(f_0)$, compared with the actual value of the error ${\rm Error}^{[2]}$ and the value of the integral $I_{\omega}^{[2]} (f_0)$ for some
			values of $n$, $\omega$ in the case of $d\sigma_n^{[2]}$}\label{figura0v}
	\end{center}
\end{table}

\begin{table}[!htp]
	\begin{center}{\footnotesize
			\begin{tabular}{llllll}\hline
				$n,\omega$ & $r^{[2]}_1(f_1)$ & $r^{[2]}_2(f_1)$
				& $r^{[2]}_3(f_1)$ & ${\rm Error}^{[2]}$ & $I_{\omega}^{[2]} (f_1)$ \\
				\hline
				$5, 0.2$ & $1.091(-14) $ &  $1.092(-14) $ &$1.090(-14)$  & $ 1.477(-15)$ &
				$4.1285...(-3)$ \\
				$10, 0.2$ & $1.212(-31) $ & $1.212(-31) $ &$1.211(-31)$  & $1.100(-32)$ &
				$4.0317...(-6)$ \\
				$15, 0.2 $ & $ 3.241(-49)$ & $ 3.243(-49)$ &$3.239(-49)$  & $ 2.327(-50)$ &
				$3.9372...(-9)$ \\
				$20, 0.2 $ & $ 3.703(-67)$ & $ 3.705(-67)$ &$3.701(-67)$  & $ 2.250(-68)$ &
				$3.8449...(-12)$\\
				\hline
				$5, 1$ & $6.844(-8) $ & $6.995(-8) $ &$6.668(-8)$  & $ 9.110(-9)$ &
				$3.3409...(-3)$ \\
				$10, 1$ & $6.217(-18) $ & $6.312(-18) $ &$6.111(-18)$  & $5.579(-19)$ &
				$3.3626...(-6)$ \\
				$15, 1 $ & $ 1.406(-28)$ & $ 1.423(-28)$ & $1.386(-28)$  & $ 9.984(-30)$&
				$3.1861...(-9)$ \\
				$20, 1 $ & $ 1.379(-39)$ & $ 1.394(-39)$ &$1.362(-39)$  & $ 8.296(-41)$&
				$3.1115...(-12)$
				\\
				\hline
				$5, 5$ & $2.900(-3) $ & $2.800(-3)$ &$1.800(-3)$  & $ 2.135(-4)$ &
				$1.1490...(-3)$ \\
				$10, 5$ & $7.964(-8) $ & $7.683(-8) $ &$5.447(-8)$  & $4.578(-9)$ &
				$1.4971...(-6)$ \\
				$15, 5 $ & $ 9.769(-13)$ & $ 9.771(-13)$ &$7.159(-13)$  & $ 4.834(-14)$ &
				$1.4612...(-9)$ \\
				$20, 5 $ & $ 7.392(-18)$ & $7.562(-18)$ &$5.702(-18)$  & $ 3.102(-19)$ &
				$1.4270...(-12)$ \\
				\hline
		\end{tabular}}
		\vskip 0.2in
		\caption{The values of the derived bounds $r_1^{[2]}(f_1),r_2^{[2]}(f_1),r_3^{[2]}(f_1)$, compared with the actual value of the error ${\rm Error}^{[2]}$ and the value of the integral $I_{\omega}^{[2]} (f_1)$ for some
			values of $n$, $\omega$ in the case of $d\sigma_n^{[2]}$}\label{figura1v}
	\end{center}
\end{table}

\begin{table}[!htp]
	\begin{center}{\footnotesize
			\begin{tabular}{llllll}\hline
				$n,\omega$ & $r^{[3]}_1(f_0)$ & $r^{[3]}_2(f_0)$
				& $r^{[3]}_3(f_0)$ & ${\rm Error}^{[3]}$ & $I_{\omega}^{[3]} (f_0)$ \\
				\hline
				$5, 0.1$ & $1.557(-12) $ & $1.496(-12)$ &$1.497(-12)$  & $ 5.250(-13)$ &
				$3.2272...(-3)$ \\
				$10, 0.1$ & $6.799(-28) $ & $6.617(-28)$  &$6.619(-28)$ & $ 1.656(-28)$ &
				$3.1516...(-6)$ \\
				$15, 0.1$ & $2.186(-44) $ & $2.139(-44)$ &$2.139(-44)$  & $4.381(-45)$ &
				$3.0777...(-9)$ \\
				$20, 0.1 $ & $ 1.288(-61)$ & $1.264(-61)$ &$1.264(-61)$  & $ 2.246(-62)$&
				$3.0056...(-12)$ \\
				\hline
				$5, 1$ & $2.779(-7) $ & $2.376(-7)$ &$2.391(-7)$ & $ 8.319(-8)$ &
				$5.3793...(-3)$ \\
				$10, 1$ & $1.057(-17) $ & $1.044(-17)$ &$1.048(-17)$  & $2.611(-18)$ &
				$5.2523...(-6)$ \\
				$15, 1 $ & $ 3.649(-29)$ & $3.372(-29)$ &$3.379(-29)$  & $ 6.898(-30)$ &
				$5.1301...(-9)$ \\
				$20, 1 $ & $ 2.129(-41)$ & $1.991(-41)$ &$1.994(-41)$  & $ 3.533(-42)$ &
				$5.0098...(-12)$ \\
				\hline
				$5, 5$ & $1.200(-2) $ & $7.400(-3)$  &$7.600(-3)$ & $ 2.455(-3)$ &
				$1.2295...(-1)$  \\
				$10, 5$ & $1.175(-9) $ & $8.798(-10)$ &$8.928(-10)$  & $2.158(-10)$ &
				$1.2007...(-4)$ \\
				$15, 5 $ & $ 1.047(-17)$ & $8.402(-18)$ &$8.488(-18)$  & $ 1.708(-19)$ &
				$1.1726...(-7)$ \\
				$20, 5 $ & $ 1.814(-26)$ &  $1.510 (-26)$& $1.522(-26)$  & $ 2.673(-27)$ &
				$1.1451...(-10)$\\
				\hline
		\end{tabular}}
		\vskip 0.2in
		\caption{The values of the derived bounds $r_1^{[3]}(f_0),r_2^{[3]}(f_0),r_3^{[3]}(f_0)$, compared with the actual value of the error ${\rm Error}^{[3]}$ and the value of the integral $I_{\omega}^{[3]} (f_0)$ for some
			values of $n$, $\omega$ in the case of $d\sigma_n^{[3]}$}\label{figura0g}
	\end{center}
\end{table}

\begin{table}[!htp]
	\begin{center}{\footnotesize
			\begin{tabular}{llllll}\hline
				$n,\omega$ & $r^{[3]}_1(f_1)$ & $r^{[3]}_2(f_1)$
				& $r^{[3]}_3(f_1)$ & ${\rm Error}^{[3]}$ & $I_{\omega}^{[3]} (f_1)$ \\
				\hline
				$5, 0.2$ & $1.117(-14) $ &$1.090(-14)$ & $1.090(-14)$ & $ 2.951(-15)$ &
				$8.8257...(-3)$ \\
				$10, 0.2$ & $1.235(-31) $ & $1.211(-31)$ &$1.211(-31)$  & $2.200(-32)$ &
				$8.0634...(-6)$ \\
				$15, 0.2 $ & $ 3.297(-49)$ & $3.239(-49)$ & $3.240(-49)$ & $ 4.651(-50)$ &
				$7.8745...(-9)$ \\
				$20, 0.2 $ & $ 3.762(-67)$ & $3.701(-67)$ & $3.702(-67)$ & $ 4.497(-68)$ &
				$7.6899...(-12)$\\
				\hline
				$5, 1$ & $7.797-8) $ & $6.666(-8)$  & $6.707(-8)$ & $ 1.785(-8)$ &
				$6.6819...(-3)$ \\
				$10, 1$ & $6.896(-18) $ & $6.110(-18)$ & $6.136(-18)$ & $1.099(-18)$ &
				$6.5253...(-6)$ \\
				$15, 1 $ & $ 1.542(-28)$ & $1.386(-28)$ &  $1.391(-28)$ & $ 1.972(-29)$&
				$6.3723...(-9)$ \\
				$20, 1 $ & $ 1.502(-39)$ & $1.362(-19)$ & $1.366(-19)$  & $ 1.641(-40)$&
				$6.2230...(-12)$
				\\
				\hline
				$5, 5$ & $5.500(-3) $ &  $1.800(-3)$ &  $1.800(-3)$& $ 3.451(-4)$ &
				$3.0645...(-3)$ \\
				$10, 5$ & $1.365(-7) $ & $5.319(-8)$ & $5.605(-8)$ & $7.546(-9)$ &
				$2.9927...(-6)$ \\
				$15, 5 $ & $ 1.593(-12)$ & $7.027(-13)$ & $7.361(-13)$ & $ 8.124(-14)$ &
				$2.9925...(-9)$ \\
				$20, 5 $ & $ 1.157(-17)$ & $5.618(-18)$ & $5.856(-18)$ & $ 5.280(-19)$ &
				$2.8540...(-12)$ \\
				\hline
		\end{tabular}}
		\vskip 0.2in
		\caption{The values of the derived bounds $r_1^{[3]}(f_1),r_2^{[3]}(f_1),r_3^{[3]}(f_1)$, compared with the actual value of the error ${\rm Error}^{[3]}$ and the value of the integral $I_{\omega}^{[3]} (f_1)$ for some
			values of $n$, $\omega$ in the case of $d\sigma_n^{[3]}$}\label{figura1g}
	\end{center}
\end{table}

\newpage

\section{Proof of the main results}

\subsection{Proof of Theorem \ref{thm:maxmod}.}

For the proof of this theorem it will be useful the following
\begin{lemma}\label{lem:tech}
Let $A$ and $B$ be two real numbers different each other and $\zeta=\rho e^{i\varphi}$ and
\[
I(\varphi)=\frac{\zeta^2+A+o\left(\frac{1}{\rho}\right)}{\zeta^2+B+o\left(\frac{1}{\rho}\right)}, \quad (\rho\to\infty)\,.
\]
If $A<B$, then there exists $\rho^*$ such that for each $\rho>\rho^*$
\[
\max_{\varphi\in [0, 2\pi)}\left|I(\varphi)\right|=\left|I\left(\frac{\pi}{2}\right)\right|
\]
and if $A>B$, then there exists $\rho^{**}$ such that for each $\rho>\rho^{**}$
\[
\max_{\varphi\in [0, 2\pi)}\left|I(\varphi)\right|=\left|I\left(0\right)\right|=\left|I\left(\pi\right)\right|.
\]
\end{lemma}

\textbf{Proof} For each $X\in\mathbb{R}$ we have
\begin{eqnarray*}
	 &&\left|\zeta^2+X+o\left(\frac{1}{\rho}\right)\right|^2=\left(\rho^2\cos{2\varphi}+X+o\left(\frac{1}{\rho}\right)\right)^2+\left(\rho^2\sin{2\varphi}+o\left(\frac{1}{\rho}\right)\right)^2\\
	&=&\rho^4+2X\rho^2\cos{2\varphi}+o\left(\rho^2\right)=\rho^2\left(\rho^2+2X\cos{2\varphi}+o\left(1\right) \right) \quad (\rho\to\infty),
\end{eqnarray*}
so we have to prove that for large enough $\rho$ it holds
\[
\frac{\rho^2+2A\cos{2\varphi}+o\left(1\right)}{\rho^2+2B\cos{2\varphi}+o\left(1\right)}\leq \frac{\rho^2-2A+o\left(1\right)}{\rho^2-2B+o\left(1\right)}
\]
for each $\theta\in[0,2\pi)$ if $A<B$
and
\[
\frac{\rho^2+2A\cos{2\varphi}+o\left(1\right)}{\rho^2+2B\cos{2\varphi}+o\left(1\right)}\leq \frac{\rho^2+2A+o\left(1\right)}{\rho^2+2B+o\left(1\right)}
\]
for each $\theta\in[0,2\pi)$ if $A>B$, i.e.
\[
\begin{array}{rl}
\left(\rho^2+2A\cos{2\varphi}+o\left(1\right)\right)\left(\rho^2-2B+o\left(1\right)\right)\\[0.1in]
-\left(\rho^2-2A+o\left(1\right)\right)\left(\rho^2+2B\cos{2\varphi}+o\left(1\right)\right)\leq 0
\end{array}
\]
for each $\theta\in[0,2\pi)$ if $A<B$
and
\[
\begin{array}{rl}
\left(\rho^2+2A\cos{2\varphi}+o\left(1\right)\right)\,\left(\rho^2+2B+o\left(1\right)\right)\\[0.1in]
-\left(\rho^2+2A+o\left(1\right)\right)\,\left(\rho^2+2B\cos{2\varphi}+o\left(1\right)\right)\leq 0
\end{array}
\]
for each $\theta\in[0,2\pi)$ if $A>B$; but it is obvious since the last two expressions are of the form
\[
2(A-B)(1+\cos{2\varphi})\rho^2+o(\rho^2) \quad _{(\rho\to\infty)}
\]
and
\[
2(B-A)(1-\cos{2\varphi})\rho^2+o(\rho^2) \quad _{(\rho\to\infty)}\,,
\]
respectively.

\textbf{Proof of Theorem \ref{thm:maxmod}}.

\textbf{(1)} Since for $m>2$, $m$ even, from (\ref{JezgroAp}) we get
\[
K^{[I]}_m(z)=(-1)^{m/2}\pi\frac{\zeta^m\left(\zeta^2+1\right)+o(\zeta^m)}{\zeta^{2m+1}(\zeta^{2}-2)+o(\zeta^{2m-1})}=\frac{(-1)^{m/2}\pi}{\zeta^{m+1}}\frac{\zeta^2+1+o(1)}{\zeta^2-2+o(1)} \,,\quad \rho\to+\infty,
\]
the statement directly follow from Lemma \ref{lem:tech}, since only the modulus of $\displaystyle \frac{\zeta^2+1+o(1)}{\zeta^2-2+o(1)}$ depends on the argument of $\zeta$ in the complex plane.
For $m=2$ we have
\[
K^{[I]}_2(z)=-\frac{\pi}{\zeta^3}\frac{\zeta^4+o(\zeta^2)}{\zeta^3-3\zeta+o(\zeta)}=-\frac{\pi}{\zeta^2}\frac{\zeta^2+o(1)}{\zeta^2-3+o(1)}\,, \quad \rho\to+\infty,
\]
and the statement again follows in the same way.

For odd $m$, from (\ref{JezgroAn}) it holds
\[
K^{[I]}_m(z)=\frac{\pi(m+2)}{m\zeta^{2m+1}}\frac{\zeta^2+\frac{m}{m+2}+o(1)}{\zeta^2-\frac{2m-2}{m}+o(1)}\,, \quad \rho\to+\infty
\]
and the proof again readily follows from Lemma \ref{lem:tech}.

\textbf{(2)}  From (\ref{Jezgro2}), we get that for $n>1$,
\[
K^{[1]}_n(z)=\frac{3\pi}{2^{2n-2}}\frac{\zeta^2+o(1)}{\zeta^2-1+o(1)}, \quad \rho\to+\infty
\]
and
\[
K^{[1]}_1(z)=3\pi\frac{\zeta^2+\frac{1}{3}+o(1)}{\zeta^2+o(1)},  \quad \rho\to+\infty,
\]
so the proof directly follows again from Lemma \ref{lem:tech}.

\textbf{(3)} From (\ref{Jezgro3}) it is clear that for $n>1$,
\[
K^{[2]}_n(z)=\frac{\pi}{2^{2n-1}\zeta^{2n-1}}\frac{\zeta^2-0.5+o(1)}{\zeta^2-1+o(1)}, \quad \rho\to+\infty
\]
and
\[
K^{[2]}_1(z)=\frac{\pi}{2\zeta^{3}}\frac{\zeta^2-0.5+o(1)}{\zeta^2+o(1)},  \quad \rho\to+\infty,
\]
so a straightforward application of Lemma \ref{lem:tech} renders the proof.

\textbf{(4)} In this case we can not directly use the results of Lemma \ref{lem:tech}, because there are two consecutive natural powers of $\rho$ in the numerator. Anyway, the modulus of the kernel  (\ref{Jezgro4}) admits the asymptotic expression ($\zeta=\rho e^{i\theta}$)
\[
\frac{\pi}{4^n\rho^{3n+1}}\frac{4\rho^{4n+2}+8\rho^{4n+1}
\cos\theta+o\left(\rho^{4n+1}\right)}{\left(\frac{1}{2}\rho^2-cos2\theta+o(1)\right)\left(\frac{1}{2}\rho^{2n}+
\cos2n\theta+o(1)\right)},\quad  \rho\to+\infty,
\]
so, we have to prove
\[
\frac{\rho+2\cos\theta+o(1)}{\left(\rho^2-2\cos{2\theta}+o(1)\right)\left(\rho^{2n}+2\cos{2n\theta}+o(1)\right)}\leq \frac{\rho+2+o(1)}{\left(\rho^2-2+o(1)\right)\left(\rho^{2n}+2+o(1)\right)}, \quad \rho\to+\infty,
\]
or, what is the same,
\[
\begin{array}{rl}
\left(\rho+2\cos\theta+o(1)\right)\left(\rho^2-2+o(1)\right)\left(\rho^{2n}+2+o(1)\right)\\[0.1in]
-\left(\rho+2+o(1) \right)\left(\rho^2-2\cos{2\theta}+o(1)\right)\left(\rho^{2n}+2\cos{2n\theta}+o(1)\right)<0
\end{array}
\]
for $\rho$ large enough and whatever $\theta\neq 0$, which is obvious since the expression above is asymptotically of the form
\[
2\rho^{2n+2}\left(\cos{\theta}-1\right)+o\left(\rho^{2n+2}\right), \quad \rho\to+\infty.
\]

\subsection{Proof of Theorem \ref{thm:errors2}.}

\textbf{(1)} This is the case where $i=1$ and $n>1$. This kernel is very similar to the kernel which appears in \cite{triocene}, so we will expand it following the analogous steps.

First, we need to state a few technical lemmas.

\begin{lemma}\label{lem:tech1}
	If $z\notin [-1,1]$, then holds the following expansion
	\begin{equation}\label{razvojceb2}
	\dfrac{1}{T_n(z)}=\sum_{k=0}^{+\infty}\beta^{[1]}_{n,k} \zeta^{-2n-k},
	\end{equation}
 with
 \begin{equation}\label{uslovi}
 \zeta=\rho e^{i\theta}, \quad \rho>1, \quad z= (\zeta+\zeta^{-1})/{2}
 \end{equation}
	where
	\begin{equation}\label{bete2}
	\beta^{[1]}_{n,k}=\left\{\begin{array}{ll}
	2(-1)^j, & k=2jn, \\
	0, & \mbox{otherwise}.
	\end{array}\right.
	\end{equation}
\end{lemma}

{\bf Proof}
	We know that if $x\in \mathbb{C}$, $|x|<1$, then
	\begin{equation}
	\dfrac{1}{(1-x)^{\nu+1}}=\sum_{k=\nu}^{+\infty}{k \choose \nu}x^{k-\nu}, \ \nu=0,1,2,...\ .
	\end{equation}
	
	Using this fact and (\ref{Ceb1}), we get
	\begin{eqnarray*}
		\dfrac{1}{T_n(z}&=&\left[\dfrac{1}{2}(\zeta^n+\zeta^{-n})
		\right]^{-1}=2\zeta^{-n}\dfrac{1}{1+\zeta^{-2n}}\\
		&=& 2\sum_{j=0}^{+\infty}(-1)^j\zeta^{-n-2nj},
	\end{eqnarray*}
	which completes the proof.\hfill\hfill\qed
	
The following lemma was proved in the paper \cite{triocene}.

\begin{lemma}\label{lem:tech2}
	If $z \notin [-1,1]$, using the notation in (\ref{uslovi}), $\varrho^{[1]}_n$ can be expanded as
	\begin{equation}\label{razvojrho2}
	\varrho^{[1]}_n(z)=\dfrac{1}{2^{3n-3}}\sum_{k=0}^{+\infty}\gamma^{[1]}_{n,k}\zeta^{-n-k-1},
	\end{equation}
	where
	\begin{equation}\label{game2}
	\gamma^{[1]}_{n,k}=\left\{\begin{array}{lll}
	\displaystyle\frac{3\pi}{2}, & k=0,2,...,2n-2, \\
	2\pi, & k=2n,2n+2,...,\\
	0, & \mbox{otherwise}.
	\end{array}\right.
	\end{equation}
\end{lemma}

Then, substituting (\ref{razvojceb2}) and (\ref{razvojrho2}) in
(\ref{Jezgro2}), we obtain
\begin{equation}\label{razvojjezgro2}
K^{[1]}_n(z)=\dfrac{1}{2^{2n-2}}\sum_{k=0}^{+\infty}\omega^{[1]}_{n,k}\zeta^{-2n-k-1},
\end{equation}
where
\begin{equation}\label{konvolucija2}
\omega^{[1]}_{n,k}=\dfrac{1}{2^{2n-2}}\sum_{j=0}^k\beta^{[1]}_{n,j}\gamma^{[1]}_{n,k-j}.
\end{equation}

Now, we are in a position to state the following result about the expansion of the remainder term $\displaystyle R^{[1]}_n(f)$, which leads to the proof of Theorem \ref{thm:errors2}, part (1).

\begin{theorem}\label{thm:remainder1}
	The remainder term $R^{[1]}_n(f)$ can be represented in the form
	\begin{equation}\label{razvojostatak2}
	R^{[1]}_n(f)=\dfrac{1}{2^{2n-2}}\sum_{k=0}^{+\infty}\alpha_{2n+k}\epsilon^{[1]}_{n,k},
	\end{equation}
	where the coefficients $\epsilon^{[1]}_{n,k}$ are independent on $f$.
	Furthermore, if $f$ is an even function, then $\epsilon^{[1]}_{n,2j+1}=0$
	$(j=0,1,...).$
\end{theorem}
	
	{\bf Proof}
Now, using (\ref{cebred}) and (\ref{razvojjezgro2}) in
	(\ref{razvojostatak2}) we obtain
	\begin{eqnarray*}
		R^{[1]}_n(f)&=&\dfrac{1}{2^{2n-2}}\dfrac{1}{2\pi i}\int_{\mathcal{E}_\rho}\left(\sum_{k=0}^{\infty}{'}\alpha_k T_k(z)\sum_{k=0}^{+\infty}\omega^{[1]}_{n,k}\zeta^{-2n-k-1}\right)dz \\
		&=&\dfrac{1}{2^{2n-2}}\sum_{k=0}^{+\infty}\left(\dfrac{1}{2\pi i}\sum_{j=0}^{+\infty}{'}\alpha_j \right.
		\left.\int_{\mathcal{E}_\rho}T_j(z)\zeta^{-2n-k-1}dz\right)\omega^{[1]}_{n,k}.
	\end{eqnarray*}
	Applying \cite[Lemma 5]{Hunter}, this reduces to
	(\ref{razvojostatak2}) with
	\begin{equation}\label{veza}
	\epsilon^{[1]}_{n,0}=\dfrac{1}{4}\omega^{[1]}_{n,0}, \
	\epsilon^{[1]}_{n,1}=\dfrac{1}{4}\omega^{[1]}_{n,1}, \
	\epsilon^{[1]}_{n,k}=\dfrac{1}{4}(\omega^{[1]}_{n,k}-\omega^{[1]}_{n,k-2}), \
	k=2,3,...\ .
	\end{equation}
	When $k$ is odd, since $\omega(t)=\omega(-t)$ it follows from
	(\ref{konvolucija2}) and Lemmas \ref{lem:tech1} and \ref{lem:tech2} that $\omega^{[1]}_{n,k}=0$,
	and hence $\epsilon^{[1]}_{n,k}=0$.\hfill\hfill\qed

Now, by using (\ref{bete2}), (\ref{game2}), (\ref{konvolucija2}), we have that, if and
only if $k=2jn$, $j\in \mathbb{N}_0$,
\begin{eqnarray*}
	 \omega^{[1]}_{n,2jn}&=&\beta^{[1]}_{n,0}\gamma^{[1]}_{n,2jn}+\beta^{[1]}_{2n}\gamma^{[1]}_{n,(2j-2)n}+...+\beta^{[1]}_{n,(2j-2)n}\gamma^{[1]}_{n,2n}+\beta^{[1]}_{n,2jn}\gamma^{[1]}_{n,0},\\
	\omega^{[1]}_{n,2jn-2}&=&\beta^{[1]}_{n,0}\gamma^{[1]}_{n,2jn-2}+\beta^{[1]}_{n,2n}\gamma^{[1]}_{(2j-2)n-2}+...+\beta^{[1]}_{n,(2j-2)n}\gamma^{[1]}_{n,2n-2},
\end{eqnarray*}
which implies
\begin{eqnarray*}
	\omega^{[1]}_{n,2jn}-\omega^{[1]}_{n,2jn-2}&=&\beta^{[1]}_{n,(2j-2)n}(\gamma^{[1]}_{n,2n}-\gamma^{[1]}_{n,2n-2})+\beta^{[1]}_{n,2jn}\gamma^{[1]}_{n,0}\\
	&=&2(-1)^{j-1}j\dfrac{\pi}{2}+4(-1)^j\dfrac{3\pi}{2}=2\pi(-1)^j,
\end{eqnarray*}
i.e.
\[
\epsilon^{[1]}_{n,2jn}=(-1)^j\dfrac{\pi}{2}.
\]
Otherwise, $\epsilon^{[1]}_{n,k}=0$ for $k\neq 2jn$. Using previously obtained
results, we get
\begin{eqnarray*}
	|R^{[1]}_n(f)|&=&\dfrac{1}{2^{2n-2}}\left|\sum_{k=0}^{+\infty}\alpha_{2n+k}\epsilon^{[1]}_{n,k} \right|=\left|\sum_{k=0}^{+\infty}\alpha_{2n+2jn}\epsilon^{[1]}_{n,2jn} \right|\\
	&\leq& \dfrac{1}{2^{2n-2}}\dfrac{\pi}{\rho^{2n}}\,\|f\|_{\mathcal{E}_{\rho}}\,\sum_{k=0}^{+\infty}\dfrac{1}{\rho^{2jn}}\\
	&=&\dfrac{1}{2^{2n-2}}\dfrac{\pi}{\rho^{2n}-1}\,\|f\|_{\mathcal{E}_{\rho}}\,,
\end{eqnarray*}
and the error bound \eqref{Druga1} easily follows.

\textbf{(2)} Using the standard expansion for $\frac{1}{1-\zeta^{-2}}$, (\ref{ro3}) can be written as
\begin{eqnarray*}
\rho^{[2]}_n(z)&=&\dfrac{\pi}{2^{3n}}\left(\sum_{p=0}^{+\infty}\zeta^{-2p}\right) \left(2\zeta^{-(n+1)}-\zeta^{-(n+3)}-\zeta^{-(3n+3)} \right)\\
&=&\dfrac{\pi}{2^{3n}}\left(2\sum_{p=0}^{+\infty}\zeta^{-(n+2p+1)}-\sum_{q=0}^{+\infty}\zeta^{-(n+2q+3)}-\sum_{r=0}^{+\infty}\zeta^{-(3n+2r+3)} \right)\\
&=&\dfrac{\pi}{2^{3n}}\left(2\zeta^{-(n+1)}+\sum_{p=1}^{n}\zeta^{-(n+1+2p)}\right),
\end{eqnarray*}
and further using (\ref{Jezgro3}), we get
\begin{eqnarray*}
K^{[2]}_n(z)&=&\dfrac{\pi}{2^{3n}}\left(2\zeta^{-(n+1)}+\sum_{p=1}^{n}\zeta^{-(n+1+2p)}\right)\dfrac{2^n\zeta^{-n}}{1+\zeta^{-2n}}\\
&=&\dfrac{\pi}{2^{2n}}\left(2\zeta^{-(2n+1)}+\sum_{p=1}^{n}\zeta^{-(2n+1+2p)}\right)\left(\sum_{j=0}^{+\infty}(-1)^j\zeta^{-2jn} \right)\\
&=&\dfrac{\pi}{2^{2n}\zeta^{2n+1}}\left(2+\sum_{p=0}^{n}{\zeta^{-2p}} \right)\left(\sum_{j=0}^{+\infty}(-1)^j\zeta^{-2jn} \right)\\
&=&\dfrac{\pi}{2^{2n}\zeta^{2n+1}}\left[\left(2+{\zeta^{-2n}}\right)\left(\sum_{j=0}^{+\infty}(-1)^j\zeta^{-2jn} \right)+\left(\sum_{p=1}^{n-1}{\zeta^{-2p}} \right)\left(\sum_{j=0}^{+\infty}(-1)^j\zeta^{-2jn} \right)\right]\\
&=&\dfrac{\pi}{2^{2n}\zeta^{2n+1}}\left(\sum_{j=0}^{+\infty}2(-1)^j\zeta^{-2nj}+\sum_{j=0}^{+\infty}(-1)^{(j)}\zeta^{-2(j+1)n}+\sum_{j=0}^{+\infty}\sum_{p=1}^{n-1}(-1)^j\zeta^{-(2jn+2p)}\right)\\
&=&\dfrac{\pi}{2^{2n}\zeta^{2n+1}}\left(1+\sum_{j=0}^{+\infty}\sum_{p=0}^{n-1}(-1)^j\zeta^{-(2jn+2p)}\right).
\end{eqnarray*}

Hence, here we also have the expansion
\begin{equation}\label{razvojjezgro3}
K^{[2]}_n(z)=\dfrac{1}{2^{2n}}\sum_{k=0}^{+\infty}\omega^{[2]}_{n,k}\zeta^{-2n-k-1},
\end{equation}

where
\begin{equation}\label{omege2}
\omega^{[2]}_{n,k}=\left\{\begin{array}{lll}
\displaystyle 2, & k=0, \\
(-1)^j, & k\in\mathbb{N}, \quad k=2nj+2p, \ j\in\mathbb{N}_0, \ p\in\{0,1,...,n-1\}, \\
0, & \mbox{otherwise}.
\end{array}\right.
\end{equation}
Here "otherwise" obviously means that $k$ is odd.

Now, based on the same principle as in (\ref{veza}), we get
\begin{eqnarray*}
\epsilon^{[2]}_{n,0}&=&\dfrac{1}{4}\omega^{[2]}_{n,0}=\dfrac{1}{2}, \\
\epsilon^{[2]}_{n,2}&=&\dfrac{1}{4}\left(\omega^{[2]}_{n,2}-\omega^{[2]}_{n,0}\right)=-\frac{1}{4}, \\
\epsilon^{[2]}_{n,k}&=&\dfrac{1}{4}\left(\omega^{[2]}_{n,k}-\omega^{[2]}_{n,k-2}\right)=\left\{\begin{array}{ll}
	\displaystyle \frac{(-1)^j}{2}, & k=2jn, \quad j\in\mathbb{N}, \\
	0, & \mbox{otherwise},\quad k\geq 1.
\end{array}\right. \\
\end{eqnarray*}
Analogously to (\ref{razvojostatak2}), we get
	\begin{equation}\label{razvojostatak3}
R^{[2]}_n(f)=\dfrac{1}{2^{2n}}\sum_{k=0}^{+\infty}\alpha_{2n+k}\epsilon^{[2]}_{n,k},
\end{equation}
and then, using (\ref{cebocena}),
\begin{eqnarray*}
	|R^{[2]}_n(f)|&=&\dfrac{1}{2^{2n}}\left|\sum_{k=0}^{+\infty}\alpha_{2n+k}\epsilon^{[2]}_{n,k} \right|=\dfrac{1}{2^{2n}}\left|-\frac{1}{4}\alpha_{2n+2}+\sum_{j=0}^{+\infty}\alpha_{2n+2jn}\epsilon^{[2]}_{n,2jn} \right|\\
	&\leq& \dfrac{1}{2^{2n}}\dfrac{\pi}{\rho^{2n}}\,\|f\|_{\mathcal{E}_{\rho}}\,\left(\frac{1}{2\rho^2}+\sum_{k=0}^{+\infty}\dfrac{1}{\rho^{2jn}}\right)\\
	&=&\dfrac{\pi}{2^{2n}}\left(\dfrac{1}{\rho^{2n}-1}+\dfrac{1}{2\rho^{2n+2}}\right)\,\|f\|_{\mathcal{E}_{\rho}}\,,
\end{eqnarray*}
which yields the bound \eqref{Druga2}.

\textbf{(3)} Using again standard expansion for $\frac{1}{1-\zeta^{-2}}$, (\ref{ro4}) can be written as
\begin{eqnarray*}
	\rho^{[3]}_n(z)&=&\dfrac{\pi}{2^{3n}}\left(\sum_{p=0}^{+\infty}\zeta^{-2p}\right) \left(2\zeta^{-(n+1)}+\zeta^{-(n+2)}+\zeta^{-(3n+2)} \right)\\
	&=&\dfrac{\pi}{2^{3n}}\left(2\sum_{p=0}^{+\infty}\zeta^{-(n+2p+1)}+\sum_{q=0}^{+\infty}\zeta^{-(n+2q+2)}+\sum_{r=0}^{+\infty}\zeta^{-(3n+2r+2)} \right)\\
	&=&\dfrac{\pi}{2^{3n}}\left(2\sum_{p=0}^{+\infty}\zeta^{-(n+2p+1)}+\sum_{q=0}^{n-1}\zeta^{-(n+2q+2)}+2\sum_{r=0}^{+\infty}\zeta^{-(3n+2r+2)} \right),
\end{eqnarray*}
thus, using (\ref{Jezgro4}),
\begin{eqnarray*}
	K^{[3]}_n(z)&=&\dfrac{\pi}{2^{3n}}\left(2\sum_{p=0}^{+\infty}\zeta^{-(n+2p+1)}+\sum_{q=0}^{n-1}\zeta^{-(n+2q+2)}+2\sum_{r=0}^{+\infty}\zeta^{-(3n+2r+2)} \right)\dfrac{2^n\zeta^{-n}}{1+\zeta^{-2n}}\\
	&=&\dfrac{\pi}{2^{2n}\zeta^{2n+1}}\left(2\sum_{p=0}^{+\infty}\zeta^{-2p}+\sum_{q=0}^{n-1}\zeta^{-(2q+1)}+2\sum_{r=0}^{+\infty}\zeta^{-(2n+2r+1)} \right)\left(\sum_{j=0}^{+\infty}(-1)^j\zeta^{-2jn} \right).\\
\end{eqnarray*}

In the product of the sums $\sum_{p=0}^{+\infty}\zeta^{-2p}$ and $\sum_{j=0}^{+\infty}(-1)^j\zeta^{-2jn}$ the coefficient which multiplies $\frac{1}{\zeta^{2mn+2l}}$, where $m\in\mathbb{N}_0$ and $l\in\{0,1,...,m-1\}$,
is equal to
\[
\sum_{i=0}^m(-1)^i=\frac{1-(-1)^{m+1}}{2},
\]
while the coefficient of $\frac{1}{\zeta^{2mn+2l+1}}$, where $m\in\mathbb{N}_0$ and $l\in\{0,1,...,m-1\}$, in the product of the sums $\sum_{q=0}^{n-1}\zeta^{-(2q+2)}+2\sum_{r=0}^{+\infty}\zeta^{-(2n+2r+2)}$ and $\sum_{j=0}^{+\infty}(-1)^j\zeta^{-2jn}$, is equal to
\[
(-1)^m+2\left(\sum_{i=0}^{m-1}{(-1)^{i}} \right)=1\,.
\]
Hence, we have again the expansion
\begin{equation}\label{razvojjezgro4}
K^{[3]}_n(z)=\dfrac{1}{2^{2n}}\sum_{k=0}^{+\infty}\omega^{[3]}_{n,k}\zeta^{-2n-k-1},
\end{equation}

where
\begin{equation}\label{omege3}
\omega^{[3]}_{n,k}=\left\{\begin{array}{lll}
\displaystyle {1-(-1)^{m+1}}, & k=2mn+2l,\quad m\in\mathbb{N}_0,\quad l\in\{0,1,...,m-1\}, \\
\displaystyle 1, & k=2mn+2l+1,\quad k>0, \quad m\in\mathbb{N}_0,\quad l\in\{0,1,...,m-1\}. \\
\end{array}\right.
\end{equation}

Now, inspired, as above, by the same arguments as in (\ref{veza}), we get
\begin{eqnarray*}
	\epsilon^{[3]}_{n,0}&=&\dfrac{1}{4}\omega^{[2]}_{n,0}=\frac{1}{2}, \\
	\epsilon^{[3]}_{n,2}&=&\dfrac{1}{4}\left(\omega^{[3]}_{n,2}-\omega^{[3]}_{n,0}\right)=0, \\
	\epsilon^{[3]}_{n,k}&=&\dfrac{1}{4}\left(\omega^{[3]}_{n,k}-\omega^{[3]}_{n,k-2}\right)=\left\{\begin{array}{ll}
		\displaystyle \frac{(-1)^m}{2}, & k=2mn, \quad j\in\mathbb{N}, \\
		0, & \mbox{otherwise},\quad k\geq 1.
	\end{array}\right. \\
\end{eqnarray*}
Analogously to (\ref{razvojostatak3}), we get
\begin{equation}\label{razvojostatak4}
R^{[3]}_n(f)=\dfrac{1}{2^{2n}}\sum_{k=0}^{+\infty}\alpha_{2n+k}\epsilon^{[2]}_{n,k},
\end{equation}
and then, (\ref{cebocena}) yields
\begin{eqnarray*}
	|R^{[3]}_n(f)|&=&\dfrac{1}{2^{2n}}\left|\sum_{k=0}^{+\infty}\alpha_{2n+k}\epsilon^{[2]}_{n,k} \right|=\dfrac{1}{2^{2n}}\left|\sum_{j=0}^{+\infty}\alpha_{2n+2jn}\epsilon^{[2]}_{n,2jn} \right|\\
	&\leq& \dfrac{1}{2^{2n}}\dfrac{\pi}{\rho^{2n}}\left(\max_{z\in
		\mathcal{E}_{\rho}}|f(z)|
	\right)\left(\sum_{k=0}^{+\infty}\dfrac{1}{\rho^{2jn}}\right)\\
	&=&\dfrac{\pi}{2^{2n}}\left(\dfrac{1}{\rho^{2n}-1}\right)\,\|f\|_{\mathcal{E}_{\rho}}\,,
\end{eqnarray*}
which provides the bound \eqref{Druga3}.

\subsection{Proof of Theorem \ref{thm:errors3}.}

\textbf{(1)} We need to estimate the quantity
\[
L^{[1]}(\mathcal{E}_{\rho})=\dfrac{1}{2\pi}\oint_{\mathcal{E}_{\rho}}\left|K^{[1]}(z)\right||dz|,
\]
where from (\ref{Jezgro2}) we obtain
\[
\left|K^{[1]}_n(z)\right|=\dfrac{\pi\sqrt{\rho^{-4n}+9+6\rho^{-2n}\cos{2n\theta}}}{2^{2n-1}\rho^n\sqrt{\left(a_2-\cos{2\theta}\right)\left(a_{2n}+\cos{2n\theta}\right)}},
\]
and $a_m$ is standardly given by
\begin{equation}
a_m=\dfrac{\rho^m+\rho^{-m}}{2}, \quad m\in\mathbb{N},
\end{equation}
 since
$|dz|=(1/\sqrt{2})\cdot\sqrt{a_2-\cos{2\theta}}\,d\theta$ (see
\cite{Hunter}). Thus, we have
\begin{equation}\label{izraz2}
\begin{array}{rl}
L^{[1]}(\mathcal{E}_{\rho})=&\displaystyle\dfrac{1}{\rho^n\cdot 2^{2n}\sqrt{2}}\,{\int_0^{2\pi}\sqrt{\dfrac{\rho^{-4n}+9+6\rho^{-2n}\cos{2n\theta}}{a_{2n}+\cos{2n\theta}}}\,d\theta}\\[0.2in]
=&\displaystyle\dfrac{1}{\rho^n\cdot 2^{2n-1}\sqrt{2}}\,{\int_0^{\pi}\sqrt{\dfrac{\rho^{-4n}+9+6\rho^{-2n}\cos{2n\theta}}{a_{2n}+\cos{2n\theta}}}\,d\theta}.
\end{array}
\end{equation}
Applying the Cauchy inequality to the last expression, we obtain
\begin{equation}\label{Locena2}
\begin{array}{rl}
L^{[1]}(\mathcal{E}_{\rho})\leq&\displaystyle\dfrac{\sqrt{\pi}}{\rho^n\cdot 2^{2n-1}\sqrt{2}}\,\sqrt{\int_0^{\pi}{\dfrac{\rho^{-4n}+9+6\rho^{-2n}\cos{2n\theta}}{a_{2n}+\cos{2n\theta}}\,d\theta}}\\[0.2in]
=&\displaystyle\dfrac{\sqrt{\pi}}{\rho^n\cdot 2^{2n-1}\sqrt{2}}\,\sqrt{\left(\rho^{-4n}+9 \right)I_0+6\rho^{-2n}I_1}\\[0.2in]
=&\dfrac{\pi}{\rho^n\cdot 2^{2n-1}}\sqrt{\dfrac{7\rho^{-2n}+9\rho^{2n}}{\rho^{4n}-1}},
\end{array}
\end{equation}
since from \cite[Eq. 3.616.1, 3.616.7]{Ryzhik} we know that
\begin{equation}\label{integraliI}
\begin{array}{rl}
I_0=\displaystyle\int_0^{\pi}\dfrac{d \theta}{a_{2n}+\cos{2n\theta}}=\dfrac{2\rho^{2n}\pi}{\rho^{4n}-1},\\[0.2in]
I_1=\displaystyle\int_0^{\pi}\dfrac{\cos{2n\theta}\,d \theta}{a_{2n}+\cos{2n\theta}}=\dfrac{2\pi}{\rho^{4n}-1}.
\end{array}
\end{equation}
and finally, the upper bound \eqref{Treca1} is obtained.

\textbf{(2)} Analogously, we have that the quantity
\[
L^{[2]}(\mathcal{E}_{\rho})=\dfrac{1}{2\pi}\oint_{\mathcal{E}_{\rho}}\left|K^{[2]}(z)\right||dz|,
\]
where from (\ref{Jezgro3}) we obtain
\[
\left|K^{[2]}_n(z)\right|=\dfrac{\pi\sqrt{\rho^{-4}+4+\rho^{-(4n+4)}-4\rho^{-2}\cos{2\theta}+2\rho^{-(2n+4)}\cos{2n\theta}-4\rho^{-(2n+2)}\cos{(2n+2)\theta}}}{2^{2n+1}\rho^n\sqrt{\left(a_2-\cos{2\theta}\right)\left(a_{2n}+\cos{2n\theta}\right)}},
\]
reduces to
\begin{equation}\label{izraz3}
\begin{array}{rl}
&L^{[2]}(\mathcal{E}_{\rho})=\displaystyle\dfrac{1}{\rho^n\cdot 2^{2n+2}\sqrt{2}}\;\times \\[0.2in] &\displaystyle\int_0^{2\pi}\sqrt{\dfrac{{\rho^{-4}+4+\rho^{-(4n+4)}-4\rho^{-2}\cos{2\theta}+2\rho^{-(2n+4)}
\cos{2n\theta}-4\rho^{-(2n+2)}\cos{(2n+2)\theta}}}{a_{2n}+\cos{2n\theta}}}\,d\theta\\[0.2in]
=&\displaystyle \dfrac{1}{\rho^n\cdot 2^{2n+1}\sqrt{2}}\;\times \\[0.2in] &\displaystyle \int_0^{\pi}\sqrt{\dfrac{{\rho^{-4}+4+\rho^{-(4n+4)}-4\rho^{-2}\cos{2\theta}+2\rho^{-(2n+4)}
\cos{2n\theta}-4\rho^{-(2n+2)}\cos{(2n+2)\theta}}}{a_{2n}+\cos{2n\theta}}}\,d\theta.
\end{array}
\end{equation}
Applying the Cauchy inequality to the last expression, we obtain
\begin{equation}
\begin{array}{rl}
&L^{[2]}(\mathcal{E}_{\rho})\leq \displaystyle\dfrac{\sqrt{\pi}}{\rho^n\cdot 2^{2n+1}\sqrt{2}}\;\times \\[0.2in] &\displaystyle\sqrt{\int_0^{\pi}{\dfrac{{\rho^{-4}+4+\rho^{-(4n+4)}-4\rho^{-2}\cos{2\theta}+2\rho^{-(2n+4)}\cos{2n\theta}-4\rho^{-(2n+2)}\cos{(2n+2)\theta}}}{a_{2n}+\cos{2n\theta}}\,d\theta}}\\[0.2in]
=&\displaystyle\dfrac{\sqrt{\pi}}{\rho^n\cdot 2^{2n+1}\sqrt{2}}\,\sqrt{\left(\rho^{-4}+4+\rho^{-(4n+4)} \right)I_0-4\rho^{-2}J_1+2\rho^{-(2n+4)}I_1-4\rho^{-(2n+2)}J_{n+1}},
\end{array}
\end{equation}
where $I_0$ and $I_1$ are given by (\ref{integraliI}) and
\begin{equation}\label{integralsJ}
\begin{array}{rl}
J_1=&\displaystyle\int_0^{\pi}\dfrac{\cos{2\theta}\,d \theta}{a_{2n}+\cos{2n\theta}},\\[0.2in]
J_{n+1}=&\displaystyle\int_0^{\pi}\dfrac{\cos{(2n+2)\theta}\,d \theta}{a_{2n}+\cos{2n\theta}}.
\end{array}
\end{equation}

The integrals \eqref{integralsJ} are equal to zero for each $n>1$. Namely, we have that
\[
\dfrac{1}{a_{2n}+\cos{2n\theta}}=\dfrac{1}{a_{2n}}\dfrac{1}{1+\frac{\cos{2n\theta}}{a_{2n}}}=\dfrac{1}{a_{2n}}\sum_{k=0}^{+\infty}(-1)^k\left(\frac{\cos{2n\theta}}{a_{2n}}\right)^k,
\]
so we can write
\begin{equation}\label{razvijeniintegrali}
\begin{array}{rl}
J_1=&\displaystyle\dfrac{1}{a_{2n}}\sum_{k=0}^{+\infty}(-1)^k\int_0^{\pi}\left(\frac{\cos{2n\theta}}{a_{2n}}\right)^k\cos{2\theta}\,d \theta\\[0.2in]
=&\displaystyle\dfrac{1}{a^{k+1}_{2n}}\sum_{k=0}^{+\infty}(-1)^k\int_0^{\pi}\cos^{k}{2n\theta}\cos{2\theta}\,d \theta, \\[0.2in]
J_{n+1}=&\displaystyle\dfrac{1}{a_{2n}}\sum_{k=0}^{+\infty}(-1)^k\int_0^{\pi}\left(\frac{\cos{2n\theta}}{a_{2n}}\right)^k{\cos{(2n+2)\theta}}\,d \theta\\[0.2in]
=&\displaystyle\dfrac{1}{a^{k+1}_{2n}}\sum_{k=0}^{+\infty}(-1)^k\int_0^{\pi}\cos^{k}{2n\theta}{\cos{(2n+2)\theta}}\,d \theta. \\[0.2in]
\end{array}
\end{equation}

Using \cite[Eq. 1.320.5, 1.320.7]{Ryzhik} we get that for each $k\in\mathbb{N}_0$ holds
\begin{equation}
\begin{array}{rl}
&\displaystyle\int_0^{\pi}\cos^{k}{2n\theta}\cos{2\theta}\,d \theta=\frac{1}{2^{k-1}}\sum_{i=0}^{[\frac{k}{2}]}{''}{k\choose i}\int_0^{\pi}\cos{\left((k-2i)\cdot 2n\theta\right)}\cos{2\theta}\,d \theta \\[0.2in]
=&\displaystyle\frac{1}{2^{k}}\sum_{i=0}^{[\frac{k}{2}]}{''}{k\choose i}\int_0^{\pi}\left(\cos{\left((2n(k-2i)+2)\theta\right)}+\cos{\left((2n(k-2i)-2)\theta\right)} \right)\,d \theta \\[0.2in]
=&\displaystyle\frac{1}{2^{k}}\sum_{i=0}^{[\frac{k}{2}]}{''}{k\choose i}\left(\dfrac{\sin{\left((2n(k-2i)+2)\theta\right)}}{2n(k-2i)+2}+\dfrac{\sin{\left((2n(k-2i)-2)\theta\right)}}{2n(k-2i)-2}\right)
\Bigg|_{0}^{\pi} \\[0.2in]
=&0,
\end{array}
\end{equation}
which together with (\ref{razvijeniintegrali}) directly implies $J_1=0$  (the ``the double prim'' in the sum denotes that the last summand has to
be halved if $k$  is even). The same happens with $J_{n+1}$, and it is obvious that the fact $n>1$ is important because it guarantees that the argument of the cosine function in the last integrals can not be equal to zero, since $2n$ can not be a divisor of 2.

Hence, we get
\begin{equation}\label{Locena3}
\begin{array}{rl}
L^{[2]}(\mathcal{E}_{\rho})\leq&\displaystyle\dfrac{\pi}{\rho^n\cdot 2^{2n+1}}\sqrt{\dfrac{\rho^{2n-4}+4\rho^{2n}+3\rho^{-2n-4}}{\rho^{4n}-1}}\,,
\end{array}
\end{equation}
and, as a consequence, the upper bound \eqref{Treca2} is reached.

\textbf{(3)} Now, the quantity
\[
L^{[3]}(\mathcal{E}_{\rho})=\dfrac{1}{2\pi}\oint_{\mathcal{E}_{\rho}}\left|K^{[3]}(z)\right||dz|,
\]
where (\ref{Jezgro4}) yields
\[
\left|K^{[3]}_n(z)\right|=\dfrac{\pi\sqrt{\rho^{-2}+4+\rho^{-(4n+2)}+4\rho^{-1}\cos{\theta}+2\rho^{-(2n+2)}\cos{2n\theta}+4\rho^{-(2n+1)}\cos{(2n+1)\theta}}}{2^{2n+1}\rho^n\sqrt{\left(a_2-\cos{2\theta}\right)\left(a_{2n}+\cos{2n\theta}\right)}},
\]
reduces to
\begin{equation}\label{izraz4}
\begin{array}{rl}
&L^{[3]}(\mathcal{E}_{\rho})=\displaystyle\dfrac{1}{\rho^n\cdot 2^{2n+2}\sqrt{2}}\;\times \\[0.2in] &\displaystyle{\int_0^{2\pi}\sqrt{\dfrac{\rho^{-2}+4+\rho^{-(4n+2)}+4\rho^{-1}\cos{\theta}+2\rho^{-(2n+2)}\cos{2n\theta}+4\rho^{-(2n+1)}\cos{(2n+1)\theta}}{a_{2n}+\cos{2n\theta}}}\,d\theta}\\[0.2in]
=&\displaystyle\dfrac{1}{\rho^n\cdot 2^{2n+1}\sqrt{2}}\;\times \\[0.2in] & \displaystyle{\int_0^{\pi}\sqrt{\dfrac{\rho^{-2}+4+\rho^{-(4n+2)}+4\rho^{-1}\cos{\theta}+2\rho^{-(2n+2)}\cos{2n\theta}+4\rho^{-(2n+1)}\cos{(2n+1)\theta}}{a_{2n}+\cos{2n\theta}}}\,d\theta}.
\end{array}
\end{equation}
Applying, as above, the Cauchy inequality to the last expression, we obtain
\begin{equation}
\begin{array}{rl}
&L^{[3]}(\mathcal{E}_{\rho})\leq\displaystyle\dfrac{\sqrt{\pi}}{\rho^n\cdot 2^{2n+1}\sqrt{2}}\;\times \\[0.2in] &\sqrt{\displaystyle\int_0^{\pi}{\dfrac{\rho^{-2}+4+\rho^{-(4n+2)}+4\rho^{-1}\cos{\theta}+2\rho^{-(2n+2)}\cos{2n\theta}
+4\rho^{-(2n+1)}\cos{(2n+1)\theta}}{a_{2n}+\cos{2n\theta}}\,d\theta}}\\[0.2in]
=&\displaystyle\dfrac{\sqrt{\pi}}{\rho^n\cdot 2^{2n+1}\sqrt{2}}\,\sqrt{\left(\rho^{-2}+4+\rho^{-(4n+2)} \right)I_0+4\rho^{-1}K_1+2\rho^{-(2n+2)}I_1+4\rho^{-(2n+2)}K_{n+1/2}},
\end{array}
\end{equation}
where $I_0$ and $I_1$ are given by (\ref{integraliI}) and
\begin{equation}\label{integralsK}
\begin{array}{rl}
K_1=&\displaystyle\int_0^{\pi}\dfrac{\cos{\theta}\,d \theta}{a_{2n}+\cos{2n\theta}},\\[0.2in]
K_{n+1/2}=&\displaystyle\int_0^{\pi}\dfrac{\cos{(2n+1)\theta}\,d \theta}{a_{2n}+\cos{2n\theta}}.
\end{array}
\end{equation}

Standard symmetry arguments show that both integrals \eqref{integralsK} are equal to zero.

%since we integrate from $0$ to $\pi$ a function whose graph is odd-symmetric with respect to the coordinate $\theta=\frac{\pi}{2}$, i.e.
%\[
%\dfrac{\cos{\theta}}{a_{2n}+\cos{2n\theta}}=-\dfrac{\cos{\left(\pi-\theta\right)}}{a_{2n}+\cos\left\{2n\left(\pi-\theta\right)\right\}}
%\]
%and
%\[
%\dfrac{\cos{(2n+1)\theta}}{a_{2n}+\cos{2n\theta}}=-\dfrac{\cos\left\{(2n+1)\left(\pi-\theta\right)\right\}}{a_{2n}+\cos\left\{2n\left(\pi-\theta\right)\right\}}
%\]
%for each $\theta\in[0,\pi]$.

Hence, we get
\begin{equation}\label{Locena}
\begin{array}{rl}
L^{[3]}(\mathcal{E}_{\rho})\leq&\displaystyle\dfrac{\pi}{\rho^n\cdot 2^{2n+1}}\sqrt{\dfrac{\rho^{2n-2}+4\rho^{2n}+3\rho^{-2n-2}}{\rho^{4n}-1}}\,,
\end{array}
\end{equation}
which yields \eqref{Treca3}.

\section*{Acknowledgements} The research of R. Orive is supported in part by the Research Project of Ministerio de Ciencia e Innovaci\'on (Spain) under grant MTM2015-71352-P.
The research of A.V. Pej\v cev and M.M. Spalevi\' c is supported in part by the Serbian
Ministry of Education, Science and Technological Development (Research Project: ``Methods
of numerical and nonlinear analysis with applications'' (\#174002)).

\newpage
\noindent
Ram\'{o}n Orive, rorive@ull.es\\
 Departamento de An\'alisis Matem\'atico, Universidad de La Laguna, Spain.
\\[\baselineskip]
Aleksandar V. Pej\v cev, apejcev@mas.bg.ac.rs\\
Department of Mathematics, University of Beograd, Faculty of Mechanical Engineering,
Kraljice Marije 16, 11120 Belgrade 35, Serbia.
\\[\baselineskip]
Miodrag M. Spalevi\'c, mspalevic@mas.bg.ac.rs\\
Department of Mathematics, University of Beograd, Faculty of Mechanical Engineering,
Kraljice Marije 16, 11120 Belgrade 35, Serbia.

\end{document}